\documentclass[a4paper,11pt]{scrartcl}
\usepackage{geometry,amsmath,amssymb,enumerate,bbm,ntheorem}
\usepackage[T1]{fontenc}
\usepackage[utf8]{inputenc}
\usepackage[all]{xy}

\usepackage[french]{babel}
\newcommand{\dueto}[1]{\textup{\textbf{(#1) }}}
\newcommand{\nequiv}{\not\equiv}
\newcommand{\subindex}[2]{\index{#1!#2}}
\newcommand{\subsubindex}[3]{\index{#1!#2!#3}}
\newcommand{\tmem}[1]{{\em #1\/}}

\newcommand{\tmop}[1]{\ensuremath{\operatorname{#1}}}

\newenvironment{enumeratenumeric}{\begin{enumerate}[1.] }{\end{enumerate}}
\newenvironment{enumerateroman}{\begin{enumerate}[i.] }{\end{enumerate}}
\newenvironment{itemizeminus}{\begin{itemize} }{\end{itemize}}
\newenvironment{proof}{\noindent\textbf{D\'emonstration\ }}{\hspace*{\fill}$\Box$\medskip}
\newtheorem{conjecture}{Conjecture}[section]
\newtheorem{corollary}[conjecture]{Corollaire}
\newtheorem{definition}[conjecture]{D\'efinition}
{\theorembodyfont{\rmfamily}}
\newtheorem{lemma}[conjecture]{Lemme}
\newtheorem{proposition}[conjecture]{Proposition}
{\theorembodyfont{\rmfamily}\newtheorem{remark}[conjecture]{Remarque}}
\newtheorem{theorem}[conjecture]{Th\'eor\`eme}

\title{Non-annulation effective et positivité locale des fibrés en droites amples adjoints}
\date{12 octobre 2007}
\author{Amaël BROUSTET}

\begin{document}

\maketitle
{\let\thefootnote\relax
\footnote{
\noindent
\textbf{A.M.S.~classification :} 14E05, 14E30, 14J30, 14J40, 14J45. 
 }}

\begin{abstract}
We show how to use effective non-vanishing to prove that Seshadri constants of some ample divisors are bigger than $1$ on smooth threefolds whose anticanonical bundle is nef or on Fano varieties of small coindice. We prove the non-vanishing conjecture of Kawamata in dimension $3$ in the case of line bundles of ``high'' volume.
\end{abstract}

\section{Introduction}
Introduites par Demailly dans \cite{Demailly1992}, les constantes de Seshadri mesurent en un point la positivité d'un fibré en droites. Si $X$ est une variété projective complexe, $L$ un fibré en droites nef et $x\in X$ un point lisse, on définit la constante de Seshadri de $L$ en $x$ par
$$\varepsilon (X,L;x) = \inf_{x\in C \subset X} \frac{L \cdot C}{\tmop{mult}_x C}.$$
Leur introduction fut à l'origine motivée par des applications en vue de la conjecture de Fujita. Elles jouissent de nombreuses propriétés, dont on trouvera un exposé au chapitre 5 du livre de Lazarsfeld (\cite{Lazarsfeld2004a}).
Bien que pouvant être arbitrairement petites, il est conjecturé que les constantes de Seshadri d'un diviseur gros et nef sur une variété projective sont minorées par $1$ pour tout point en position très générale, c'est-à-dire en dehors d'une union dénombrable de sous-variétés strictes.
\begin{conjecture}
  {\dueto{{\cite{Lazarsfeld2004a}}, conjecture 5.2.4}}Soit $X$ une vari\'et\'e
  projective, $L$ un diviseur gros et nef sur $X$. Alors
  \[ \varepsilon (X, L ; x) \geqslant 1 \]
  pour tout point $x \in X$ en dehors d'une union d\'enombrable de
  sous-vari\'et\'es strictes de $X$.
\end{conjecture}
\`A l'heure actuelle, on dispose de la minoration suivante :
\begin{theorem}
 {\dueto{Ein, Küchle, Lazarsfeld, {\cite{Ein1995}}}}Soit $X$ une
  vari\'et\'e projective de dimension $n$ et $L$ un diviseur gros et nef sur $X$.
  Pour tout point en position tr\`es g\'en\'erale $x \in X$, on a
  \[ \varepsilon (X, L ; x) \geqslant \frac{1}{n} . \]
\end{theorem}
Dans le cas des fibrés en droites amples, il existe une amélioration due à Nakamaye \cite{Nakamaye2005} de l'ordre de $\frac{1}{n^2}$. Dans le cas particulier de la dimension $3$, l'amélioration est plus substantielle.

\begin{theorem}
  {\dueto{Nakamaye, {\cite{Nakamaye2005}} theorem 0.2}}Soit $X$ une
  vari\'et\'e projective de dimension $3$ et $A$ un diviseur ample sur $X$.
  Pour tout point en position tr\`es g\'en\'erale $x \in X$, on a
  \[ \varepsilon (X, A ; x) \geqslant \frac{1}{2} . \]
\end{theorem}

On prouve dans la suite de cet article une minoration optimale des constantes de Seshadri pour ``beaucoup'' de diviseurs amples en dimension $3$.
\begin{theorem}\label{Th}
Soit $X$ une variété projective de dimension $3$.
 \begin{enumerate}
  \item Supposons $X$ lisse et son diviseur anticanonique nef. Soit $L$ un diviseur ample sur $X$. Pour tout point $x$ en position très générale sur $X$, la constante de Seshadri de $L$ en $x$ vérifie
$$\varepsilon(X,L;x) \geqslant 1.$$
\item Supposons $X$ à singularités canoniques et son diviseur canonique nef. Soit $L$ un diviseur gros et nef sur $X$. Pour tout point $x$ en position très générale sur $X$, la constante de Seshadri de $K_X + L$ en $x$ vérifie
$$\varepsilon(X,K_X + L;x) \geqslant 1.$$
\item Supposons $X$ lisse, $L$ un diviseur ample tel que le diviseur
$$ A = L - K_X$$
soit ample et vérifie $\sqrt[3]{A^3}>3$.

On a alors pour tout point $x$ en position très générale sur $X$
$$\varepsilon(X,L;x) \geqslant 1.$$
Si de plus $X$ n'est pas uniréglée, on peut se contenter de supposer que $\sqrt[3]{A^3} \geqslant 3$.
 \end{enumerate}
\end{theorem}

En dimension supérieure, le même type de technique permet d'obtenir un résultat similaire pour les variétés de Fano d'indice grand.
 \begin{theorem}
  \label{coind-2}

   Soit $X$ une vari\'et\'e projective de
   dimension au moins $3$, presque de Fano, factorielle, \`a singularit\'es
   terminales et $L$ un diviseur ample sur $X$. Supposons que
   \[ - K_X \equiv (n - c + 1) H \]
   pour un diviseur entier $H$ et un entier $c \leqslant 3$. Pour tout point $x$ en position tr\`es
   g\'en\'erale, on a
   \[ \varepsilon (X, L ; x) \geqslant 1. \]
\end{theorem}

Il est aussi possible de minorer les constantes de Seshadri du diviseur fondamental des variétés de Fano de dimension $4$.
\begin{theorem}
  \label{fano4}Soit $X$ une vari\'et\'e de Fano lisse de dimension $4$, $H$ un
  diviseur ample v\'erifiant $- K_X \sim rH$ pour un entier $r$ alors pour un
  point $x$ en position tr\`es g\'en\'erale, on a
  \begin{eqnarray*}
    \varepsilon (X, H ; x) \geqslant 1. &  & 
  \end{eqnarray*}
\end{theorem}

Les constantes de Seshadri mesurent une propriété de positivité forte mais asymptotique d'un fibré en droites. Les résultats précédents sont obtenus à l'aide d'une propriété de positivité plus faible mais non asymptotique : l'existence d'une section globale non nulle d'un fibré en droites. En dimension $3$ on est ainsi à même de se ramener au cas des surfaces (voir le lemme \ref{Lem}).

Une telle propriété de non-annulation ne peut être espérée pour un fibré en droites ample quelconque, cependant Kawamata a émis la conjecture suivante pour un fibré en droites adjoint. On peut à ce sujet remarquer que les diviseurs amples apparaissant aux théorèmes \ref{Th} et \ref{coind-2} sont tous adjoints.
\begin{conjecture}
  {\dueto{{\cite{Kawamata2000}}, conjecture 2.1}}Soit
  $X$ une vari\'et\'e projective normale, $B$ un $\mathbbm{R}$-diviseur
  effectif sur $X$ tel que $(X, B)$ soit une paire klt et $D$ un diviseur de
  Cartier sur $X$. Supposons que $D$ soit nef et que $H = D - (K_X + B)$ soit
  nef et gros. Alors le fibr\'e en droites $\mathcal{O}_X (D)$ a une section
  globale non nulle : $H^0 (X, \mathcal{O}_X (D)) \neq 0$.
\end{conjecture}

Dans le cas des courbes, cette conjecture est une simple application du théorème de Riemann-Roch. Le cas de la dimension $2$ a été prouvé par Kawamata (\cite{Kawamata2000}).

En dimension $3$, seuls des résultats partiels ont été obtenus. Il est bien connu que dans le cas des variétés à singularités terminales dont le diviseur anticanonique est nef, la pseudo-effectivité de la seconde classe de Chern suffit à établir la conjecture. Comme l'a noté Xie dans \cite{Xie2004}, différents auteurs ont traité de ce sujet (voir théorème \ref{c2}).
Pour une variété lisse de dimension $3$, on prouve dans la suite que les diviseurs amples adjoints à un $\mathbf{Q}$-diviseur ample de grand volume vérifient la conjecture de Kawamata.

\begin{theorem}\label{non-annulation}\label{nonannulationdim3}

  Soit $X$ une vari\'et\'e projective lisse de
  dimension $3$, $L$ un diviseur ample sur $X$. Supposons que
  $\sqrt[3]{L^3} > 3$ et que $K_X + L$ soit nef. Alors $H^0 (X,
  \mathcal{O}_X^{} (K_X + L)) \neq 0$.

\end{theorem}

On peut noter que la conjecture de Kawamata est une version plus forte d'une conjecture que Beltrametti et Sommese avaient précédemment énoncée.

\begin{conjecture}
  {\dueto{{\cite{Beltrametti1995}}, conjecture 7.2.7}}Soit
  $X$ une vari\'et\'e projective lisse de dimension $n$,  $L$ un diviseur ample sur $X$. Supposons que le diviseur $D=K_X + (n-1)L$ soit nef. Alors le fibr\'e en droites $\mathcal{O}_X (D)$ a une section
  globale non nulle : $H^0 (X, \mathcal{O}_X (D)) \neq 0$.
\end{conjecture}

Un corollaire immédiat du théorème \ref{non-annulation} est le résultat suivant.

\begin{corollary}\label{non-annulation-Beltrammeti}
 Soit $X$ une variété projective lisse de dimension $3$ et $L$ un diviseur ample sur $X$ vérifiant $\sqrt[3]{L^3}>\frac{3}{2}$.
 Alors le fibré en droites $\mathcal{O}_X(K_X + 2L)$ a une section globale non nulle : $$H^0 (X, \mathcal{O}_X (K_X + 2L)) \neq 0.$$
\end{corollary}

Mentionnons que la conjecture de non-annulation de Beltrametti et Sommese a été récemment prouvée en dimension $3$ par Fukuma dans \cite{Fukuma2006}.

Pour les variétés de dimension supérieure à $3$, la conjecture de Kawamata ne permet pas de se ramener à une situation où il existe une minoration optimale des constantes de Seshadri (le cas les surfaces). Il existe cependant pour les variétés de Fano une conjecture plus ancienne, la conjecture des éléphants de Reid. Cette conjecture postule l'existence de membres effectifs peu singuliers dans le système linéaire du diviseur fondamental de la variété de Fano, c'est-à-dire du plus petit diviseur ample entier dont le diviseur anticanonique est un multiple. Pour les variétés dont le diviseur anticanonique est un multiple suffisamment grand de ce diviseur fondamental, on peut ainsi minorer les constantes de Seshadri du diviseur fondamental par récurrence sur la dimension, en se restreignant à ce membre peu singulier (proposition \ref{anticoind}).

En dimension $3$ comme en dimension supérieure, les constantes de Seshadri de certains diviseurs amples ne peuvent pas être minorées de façon optimale en utilisant un résultat de non-annulation. Pour un tel diviseur ample $A$, on étudie les arêtes du cône de Mori $(K_X + A)$-négatives et on utilise les résultats de classifications connus dans ces cas.

\paragraph{Remerciements.} Je remercie vivement Laurent Bonavero pour l'aide et le soutien qu'il m'a apportés tout au long de ce travail.

\section{Le cadre}
Dans la suite de cet article, toutes les vari\'et\'es consid\'er\'ees seront
suppos\'ees complexes. Si le contraire n'est pas sp\'ecifi\'e, elles seront de
plus suppos\'ees projectives.

\subsection{Diviseurs et systèmes linéaires}

Un diviseur de Weil sur une variété normale est une somme formelle finie à coefficients
entiers d'hypersurfaces. Un diviseur de Cartier est une section globale du faisceau
$\mathcal{K}_X^*/\mathcal{O}_X^*$.
On notera ``$\equiv$'' l'\'equivalence num\'erique entre diviseurs de Cartier,
``$\sim$'' l'\'equivalence lin\'eaire. Une somme formelle finie
d'hypersurfaces $D = \sum a_i D_i$ dont les coefficients sont des rationnels
est appel\'ee un $\mathbbm{Q}$-diviseur de Weil. Si un multiple non nul de ce diviseur
est un diviseur de Cartier, on dira que c'est un $\mathbbm{Q}$-diviseur de
Cartier. Lorsque l'on voudra insister sur le fait qu'un diviseur de Weil est
\`a coefficients entiers, on dira qu'il est entier.

On note ``$\sim_{\mathbbm{Q}}$'' la $\mathbbm{Q}$-\'equivalence lin\'eaire :
deux $\mathbbm{Q}$-diviseurs de Cartier $D_1$ et $D_2$ sont lin\'eairement
\'equivalents s'il existe un entier $m > 0$ tel que $mD_1$ et $mD_2$ soient des
diviseurs de Cartier entiers v\'erifiant $mD_1 \sim mD_2$. On peut d\'efinir
l'intersection d'un $\mathbbm{Q}$-diviseur de Cartier $D$ avec une courbe en
posant $D \cdot C = \frac{1}{m} (mD) \cdot C$ pour un entier $m$ tel que $mD$
soit un diviseur de Cartier entier. \

Un $\mathbbm{Q}$-diviseur de Cartier $D$ ou un fibr\'e en droites sur une
vari\'et\'e $X$ est dit {\subindex{fibr\'e en
droites}{nef}}{\tmem{nef}}{\subindex{diviseur}{nef}} si son intersection avec
toute courbe est positive, {\subindex{fibr\'e en
droites}{gros}}{\tmem{gros}}{\subindex{diviseur}{gros}} si $h^0 (X,
\mathcal{O}_X (kD)) \sim \alpha k^{\dim X}$ pour un r\'eel $\alpha > 0$ et
{\subindex{fibr\'e en
droites}{pseudo-effectif}}{\tmem{pseudo-effectif}}{\subindex{diviseur}{pseudo-effectif}}
si sa classe num\'erique appartient \`a l'adh\'erence du c\^one des diviseurs
gros dans $N^1 (X)_{\mathbbm{R}}$ (défini plus bas). On d\'efinit le {\index{volume}}volume d'un
$\mathbbm{Q}$-diviseur de Cartier $D$ par
\[ \tmop{vol} (D) = \limsup_{m \gg 0} \frac{h^0 (X, \mathcal{O}_X
   (kmD))}{(km)^{\dim X} / (\dim X) !}, \]
pour un entier $k>0$ tel que $kD$ soit Cartier.

Le lieu de base ensembliste d'une série linéaire $V$ sera noté $\tmop{Bs}(V)$, l'idéal de base de cette même série linéaire sera noté $\mathfrak{b}(V)$.

Un diviseur de Cartier $D$ sera dit {\subindex{fibr\'e en
droites}{adjoint}}adjoint{\subindex{diviseur}{adjoint}} s'il est de la forme
$D = K_X + A$ pour un $\mathbbm{Q}$-diviseur de Cartier ample $A$. On d\'efinit de m\^eme un fibr\'e en
droites adjoint.

Une propriété ``locale'' et l'expression ``localement'' feront toujours référence à la topologie de Zariski.

  On d\'efinit la relation d'\'equivalence num\'erique pour les courbes, qu'on
  notera aussi ``$\equiv$'' par $C_1 \equiv C_2$ si et seulement si $D \cdot
  C_1 = D \cdot C_2$ pour tout diviseur de Cartier $D$ sur $X$.
  
  On note $N^1 (X)_{\mathbbm{Z}}$ le quotient par la relation d'\'equivalence
  num\'erique du $\mathbbm{Z}$-module libre engendr\'e par les diviseurs de
  Cartier. On note $N^1 (X)_{\mathbbm{Q}} = N^1 (X)_{\mathbbm{Z}}
  \otimes_{\mathbbm{Z}} \mathbbm{Q}$, et $N^1 (X)_{\mathbbm{R}} = N^1
  (X)_{\mathbbm{Z}} \otimes_{\mathbbm{Z}} \mathbbm{R}$. On construit de
  mani\`ere similaire $N_1 (X)_{\mathbbm{R}}$ et $N_1 (X)_{\mathbbm{Q}}$, les
  produits tensoriels par $\mathbbm{R}$ et $\mathbbm{Q}$ du quotient par la
  relation d'\'equivalence num\'erique du $\mathbbm{Z}$-module libre
  engendr\'e par les courbes irréductibles de $X$. Le c\^one convexe engendr\'e dans $N_1
  (X)_{\mathbbm{R}}$ par les classes de courbes effectives sera not\'e $\tmop{NE} (X)$.
  Son adh\'erence sera not\'ee $\overline{\tmop{NE}} (X)$. Si $D$ est un
  $\mathbbm{Q}$-diviseur de Cartier sur $X$, on note $\overline{\tmop{NE}}
  (X)_{D \geqslant 0}$ (resp. $\overline{\tmop{NE}} (X)_{D \leqslant 0}$) la
  partie du c\^one s'intersectant positivement (resp. n\'egativement) avec le
  $\mathbbm{Q}$-diviseur $D$.

Un invariant important d'une ar\^ete $K_X$-n\'egative $R$ du cône de Mori $\overline{\tmop{NE}} (X)$ est sa \emph{longueur}, c'est-à-dire 
 l'entier
   \[ l (R) = \min \{- K_X \cdot C | [C] \in R \text{ et } C \text{ est
      rationnelle} \}. \]

On aura besoin par la suite de la notion de paire.

\subsection{Langage des paires}

Une paire{\index{paire}}\index{frontière}\label{paire} est un couple $(X, \Delta)$ o\`u $X$
  est une vari\'et\'e projective normale de dimension $n$ et $\Delta$ un
  $\mathbbm{Q}$-diviseur de Weil dont les coefficients sont tous compris, au
  sens large, entre $0$ et $1$. On dit que le diviseur $\Delta$ est une
  fronti\`ere. On requiert de plus que $K_X + \Delta$ soit
  $\mathbbm{Q}$-Cartier, le diviseur canonique $K_X$ \'etant n'importe quel
  diviseur de Weil dont la restriction \`a la partie lisse de $X$ est le
  diviseur d'une $n$-forme r\'eguli\`ere.
  
  Une log-r\'esolution{\index{log-r\'esolution}}
  d'une paire $(X, \Delta)$ est un morphisme propre, birationnel $f : Y
  \rightarrow X$ tel que $Y$ soit lisse, $\tmop{Exc} (f)$ soit un diviseur et
  $\tmop{Exc} (f) \cup \tmop{Supp} (f^{- 1}_{\ast} \Delta)$ soit un diviseur
  \`a croisement normaux simples.

Pour une log-r\'esolution $f : Y \rightarrow X$ d'une paire $(X, \Delta)$, on
d\'efinit le $\mathbbm{Q}$-diviseur de discr\'epance $\sum a (X, \Delta, E) E$
comme le diviseur v\'erifiant $K_Y = f^{\ast} (K_X + \Delta) + \sum a (X,
\Delta, E) E$ o\`u la somme porte sur tous les diviseurs irr\'eductibles $E
\subset Y$. Pour rendre ce diviseur unique, on suppose qu'un diviseur $E$ non-exceptionnel
vérifie $a(X, \Delta, E) \neq 0$ si et seulement il existe un diviseur premier $D$ sur $X$,
de coefficient $d\neq 0$ dans $\Delta$ vérifiant $E = f^{-1}_* D$. Dans ce cas on demande de plus
que $a(X, \Delta, E) = -d$.

Les singularit\'es de la paire $(X, \Delta)_{}$ peuvent \^etre class\'ees selon
les coefficients $a (X, \Delta, E)$, appel\'es
discr\'epances{\subindex{discr\'epance}{d'un diviseur}} du diviseur $E$. Ces dernières
ne dépendent pas de la log-résolution choisie.

\label{singpaire}
 On d\'efinit la
  discr\'epance{\subindex{discr\'epance}{d'une paire}} de $(X, \Delta)$ par :
  \[ \tmop{discrep} (X, \Delta) = \inf_E \{a (X, \Delta, E) \mid E 
  \textrm{ est un diviseur exceptionnel au-dessus de } X\}. \]
  On dit que la paire $(X, \Delta)$ est :
  \begin{itemizeminus}
    \item terminale {\subindex{paire}{terminale}}si $\tmop{discrep} (X,
    \Delta) > 0$ ; si $\Delta = 0$ la vari\'et\'e $X$ est dite \`a
    singularit\'es terminales{\subindex{vari\'et\'e}{\`a singularit\'es
    terminales}},
    
    \item canonique{\subindex{paire}{canonique}} si $\tmop{discrep} (X,
    \Delta) \geqslant 0$ ; si $\Delta = 0$ la vari\'et\'e $X$ est dite \`a
    singularit\'es canoniques{\subindex{vari\'et\'e}{\`a singularit\'es
    canoniques}},
    
    \item log-terminale{\subindex{paire}{log-terminale}} ou klt
    {\subindex{paire}{klt ou Kawamata log-terminale}}(pour Kawamata
    log-terminale) si $\tmop{discrep} (X, \Delta) > - 1$ et la partie
    enti\`ere de $\Delta$, $[\Delta]$, vaut $0$ ($\Delta$ est une fronti\`ere
    stricte),
    
    \item purement log-terminale{\subindex{paire}{purement log-terminale}},
    abr\'eg\'e en plt, si $\tmop{discrep} (X, \Delta) > - 1$,
    
    \item log-canonique{\subindex{paire}{log-canonique}}, abr\'eg\'e en lc, si
    $\tmop{discrep} (X, \Delta) \geqslant - 1$.
  \end{itemizeminus}

  Une paire $(X, \Delta)$ log-canonique est dite de
  Fano{\subindex{paire}{de Fano}} si le $\mathbbm{Q}$-diviseur de Cartier $-
  (K_X + \Delta)$ est ample. Si le diviseur $- (K_X + \Delta)$ est seulement
  gros et nef, on dira que $(X, \Delta)$ est presque de
  Fano{\subindex{paire}{presque de Fano}}. Une paire sera dite de
  Calabi-Yau{\subindex{paire}{de Calabi-Yau}} si $K_X + \Delta \equiv 0$ et
  de type g\'en\'eral{\subindex{paire}{de type g\'en\'eral}} si $K_X
  + \Delta$ est gros.
  
  Pour une paire presque de Fano $(X, \Delta)$, on appelle indice{\index{indice}} de la paire $(X, \Delta)$ le plus grand
rationnel $r$ tel que $- (K_X + \Delta) \equiv rH$ pour un diviseur de Cartier
entier gros et nef $H$. Un tel diviseur est appel\'e un diviseur fondamental
de la paire $(X, \Delta)$. Il est unique \`a \'equivalence num\'erique pr\`es.
Le coïndice{\index{coïndice}} de $(X, \Delta)$ est le rationnel $n + 1 - r$.

\section{Non-annulation effective}
On s'intéresse dans cette section aux résultats d'existence de sections globales non nulles pour certains fibrés en droites adjoints amples. Le point de vue adopté ici est tourné vers les applications à la positivité locale. On prouve la conjecture de Kawamata en dimension $3$ dans le cas des diviseurs gros et nefs adjoints à un diviseur ample de ``grand volume'' à la sous-section \ref{grandvolume}, les autres sous-sections consistant en quelques rappels.

\subsection{Pseudo-effectivité de la seconde classe de Chern}

Il est bien connu que la pseudo-effectivité de $c_2(X)$ implique la conjecture de non-annulation dans le cas des variétés de dimension $3$ dont le diviseur anticanonique est nef.

\begin{lemma}
  \label{secglob}Soit $X$ une vari\'et\'e projective \`a singularit\'es
  terminales de dimension $3$ dont le diviseur anticanonique est nef et dont
  la deuxi\`eme classe de Chern est pseudo-effective. Tout fibr\'e en droites
  $L$ gros et nef sur $X$ a une section globale non nulle.
\end{lemma}

\begin{proof}
  D'apr\`es le th\'eor\`eme d'annulation de Kawamata-Viehweg, le diviseur $-
  K_X$ \'etant nef, la cohomologie sup\'erieure de $\mathcal{O}_X (L)$ est
  nulle. On en d\'eduit gr\^ace au th\'eor\`eme de Riemann-Roch que $h^0 (X,
  \mathcal{O}_X (L)) = \chi (\mathcal{O}_X (L)) = \deg (\tmop{ch} (L) \cdot
  \tmop{td} (T_X))_3$, o\`u le $3$ en indice indique qu'on ne prend que la
  composante en degr\'e $3$ du produit. 
  On a alors (voir
  {\cite{Hartshorne1977}} p. 432])
  \[ h^0 (X, \mathcal{O}_X (L)) = \frac{1}{6} c_1^3 (L) + \frac{1}{4} c_1 (X)
     \cdot c_1^2 (L) + \frac{1}{12} c_1 (L) \cdot (c_1^2 (X) + c_2 (X)) +
     \chi (\mathcal{O}_X ). \]
     Or d'après \cite{Kawamata1986}, lemma 2.2 et lemma 2.3 (voir plus généralement l'ensemble de la section 2, on pourra aussi consulter \cite{Reid1987} section 8) on a $\chi (\mathcal{O}_X ) \geqslant \frac{1}{24} c_1 (X) \cdot c_2 (X)$, d'où
     \[ h^0 (X, \mathcal{O}_X (L)) \geqslant \frac{1}{6} c_1^3 (L) + \frac{1}{4} c_1 (X)
     \cdot c_1^2 (L) + \frac{1}{12} c_1 (L) \cdot (c_1^2 (X) + c_2 (X)) + \frac{1}{24} c_1 (X) \cdot c_2 (X). \]
  Puisque $c_1^3 (L) > 0$, on en d\'eduit que $h^0 (X, \mathcal{O}_X (L)) >
  0$ si les produits qui font intervenir $c_2 (X)$ sont positifs. Par
  hypoth\`ese, la deuxi\`eme classe de Chern $c_2 (X)$ est pseudo-effective et
  a donc une intersection positive ou nulle avec toute classe nef.
\end{proof}

Comme le note Xie dans {\cite{Xie2004}}, depuis Miyaoka, diff\'erents auteurs
ont trait\'e de la pseudo-effectivit\'e de la seconde classe de Chern des
vari\'et\'es dont le diviseur anticanonique est nef. Avant d'\'enoncer ces
r\'esultats (regroup\'es dans le th\'eor\`eme \ref{c2}), rappelons la notion
de dimension num\'erique d'un diviseur.

\begin{definition}
  La dimension num\'erique{\index{dimension num\'erique}} $\nu (X, D)$ d'un
  diviseur $D$ nef est le plus grand entier $k$ tel que $D^k \nequiv 0$. La
  dimension num\'erique $\nu (X)$ de $X$ est la dimension num\'erique de $-
  K_X$.
\end{definition}

\begin{theorem}
  \label{c2}Soit $X$ une vari\'et\'e projective \`a singularit\'es terminales
  de dimension $3$. Dans les cas suivants, la seconde classe de Chern de $X$
  est pseudo-effective{\subsubindex{classe de Chern}{seconde
  --}{pseudo-effective}} :
  \begin{itemizeminus}
    \item $\nu (X) = 0$ (Miyaoka, {\cite{Miyaoka1987}} theorem 1.1),
    
    \item $\nu (X) = 1$ (Keel, Matsuki et McKernan, {\cite{Keel2004}}
    corollary 6.2),
    
    \item $\nu (X) = 2$, $X$ lisse et $q (X) \neq 0$ (Xie {\cite{Xie2004}}
    theorem 3.12 et proposition 2.4),
    
    \item $\nu (X) = 3$ (Koll\'ar, Miyaoka, Mori et Takagi {\cite{Kollar2000}}
    page 5).
  \end{itemizeminus}
\end{theorem}

Pour les variétés de dimension $3$ minimales, c'est-à-dire à singularités canoniques et dont le diviseur canonique est nef, un calcul similaire à celui du lemme \ref{secglob} permet de prouver la conjecture de non-annulation (voir \cite{Kawamata2000}). Là encore, il s'agit d'utiliser le théorème de Riemann-Roch et l'inégalité de Miyaoka sur les classes de Chern de la variété.

Pour une variété de dimension $3$ quelconque, il n'est plus possible de procéder ainsi. Pour un diviseur adjoint à un diviseur de ``grand volume'', on va utiliser la sous-adjonction de Kawamata pour résoudre le problème sur une sous-variété puis étendre la section construite à la variété entière.

\subsection{Idéaux multiplicateurs et sous-adjonction}

On utilisera librement le langage des idéaux multiplicateurs tel que présenté dans \cite{Lazarsfeld2004b}.

Si on note $c = \tmop{lct} (X, D)$ le seuil log-canonique d'un
$\mathbbm{Q}$-diviseur effectif $D$ sur une vari\'et\'e lisse $X$, on peut
consid\'erer les composantes irréductibles du lieu des zéros de l'id\'eal
multiplicateur $\mathcal{J}(cD)$. Sous certaines conditions, ces
sous-vari\'et\'es sont normales et peuvent \^etre munies d'une structure de
paire qui fait apparaitre leur (log)-diviseur canonique comme la restriction
du (log)-diviseur canonique de la paire $(X, D)$. Les rappels suivants sont
\'enonc\'es dans le cadre g\'en\'eral des paires (en particulier, on ne suppose
pas que $X$ est lisse).

\subsubsection{Centres log-canoniques}

Soit $(X, D)$ une paire. On s'int\'eresse au lieu o\`u la paire
$(X, D)$ n'est pas klt. On note
\[ \tmop{nklt} (X, D) =\{x | (X, D) \text{ n'est pas klt au voisinage de }
   x\}. \]
\begin{definition}
  {\dueto{Centre de singularit\'es log-canoniques}}{\index{centre de singularités log-canoniques}}On appelle centre de singularit\'es
  log-canoniques l'image $W$ dans $X$ d'un diviseur irr\'eductible de
  discr\'epance $- 1$ sur un mod\`ele birationnel de $X$, telle que la
  paire $(X,D)$ soit log-canonique au point générique de $W$.
\end{definition}

On peut remarquer que si la paire $(X, D)$ est klt, il n'existe pas de centre
de singularit\'es log-canoniques. Plus exactement, si $(X, D)$ est log-canonique, $\tmop{nklt} (X, D)$ est
\'egal \`a l'union des centres de singularit\'es log-canoniques de $X$.

Pour une paire log-canonique, ces centres sont en nombre fini. En consid\'erant une log-r\'esolution $f : X'
\rightarrow X$ de la paire $(X, D)$, on les obtient comme l'image d'une
intersection quelconque de diviseurs de discr\'epance $- 1$.

\begin{definition}
  {\subindex{centre de singularités log-canoniques}{maximal}}{\dueto{Centre
  maximal, centre minimal}}{\subindex{centre de singularités log-canoniques}{minimal}}On d\'efinit un centre de singularit\'es
  log-canoniques maximal comme un \'el\'ement maximal pour l'inclusion. On
  d\'efinit un centre de singularit\'es log-canoniques minimal comme un
  \'el\'ement minimal pour l'inclusion.
\end{definition}

Parmi les centres minimaux, on peut consid\'erer la classe {\sl a priori} beaucoup
plus restreinte des centres de singularit\'es log-canoniques exceptionnels.

\begin{definition}
  {\dueto{Centre de singularit\'es log-canoniques
  exceptionnel}}{\subindex{centre de singularités log-canoniques}{exceptionnel}}Soit $(X, D)$ une paire log-canonique, $f : X' \rightarrow
  X$ une log-r\'esolution de la paire $(X, D)$. Un centre de singularit\'es
  log-canoniques $W$ est dit exceptionnel si les deux propri\'et\'es suivantes
  sont v\'erifi\'ees :
  \begin{enumerateroman}
    \item il existe un unique diviseur $E_W$ de discr\'epance $- 1$ sur $X'$
    dont l'image dans $X$ est $W$,
    
    \item pour tout diviseur $E' \neq E_W$ sur $X'$ de discr\'epance $- 1$, $f
    (E') \cap W = \emptyset$.
  \end{enumerateroman}
\end{definition}

On remarque qu'un centre de singularit\'es log-canoniques exceptionnel est une
composante connexe du lieu non-klt de la paire $(X, D)$. Cette derni\`ere
propri\'et\'e nous sera fort utile pour construire des sections non nulles de
fibr\'es en droites sur le lieu non-klt de certaines paires $(X, D)$.

De plus, l'image dans $X'$ de l'intersection de deux diviseurs de discrépance $-1$ étant un centre de singularités log-canonique, un centre de singularités log-canoniques minimal vérifiant la première propriété de la définition précédente est exceptionnel.

\subsubsection{Un lemme de perturbation : le ``tie-breaking''}

Les centres de singularit\'es log-canoniques exceptionnels paraissent tr\`es
sp\'eciaux parmi les centres minimaux, il n'en est en fait rien lorsque la
paire $(X, 0)$ est klt comme le montre le th\'eor\`eme suivant. Ce r\'esultat
est attribu\'e \`a Miles Reid ({\cite{Reid1983}}, 1.4) par J\'anos Koll\'ar
dans {\cite{Kollar2007}}. La version \'enonc\'ee ici l'est sous une forme
l\'eg\`erement plus forte que dans {\cite{Kollar2007}}. Cependant, la preuve
est similaire.

\begin{theorem}
  {\dueto{voir par exemple {\cite{Kollar2007}}\label{perturbation}, proposition 8.7.1}}Soit $(X, \Delta)$ une paire klt
  et $D$ un {$\mathbb{Q}$-diviseur} {$\mathbb{Q}$-Cartier} effectif tel que $(X, \Delta + D)$ soit
  log-canonique et non klt. On note $W$ un centre de singularit\'es log-canoniques
  minimal pour la paire $(X, \Delta + D)$ et $H$ un diviseur de Cartier ample sur $X$.

  Pour tout rationnel $1 \gg r > 0$, il existe des rationnels $0 \leqslant c_1
  \leqslant r$ et $0 \leqslant c_2 \leqslant r$ et un $\mathbbm{Q}$-diviseur
  effectif $A \sim_{\mathbbm{Q}} c_1 H$ tels que la paire $(X, \Delta + (1 -
  c_2) D + A)$ soit log-canonique et $W$ soit un centre de singularit\'es
  log-canononiques exceptionnel pour $(X, \Delta + (1 - c_2) D + A)$.
\end{theorem}

\begin{proof}
Si $W$ est de codimension $1$, alors $W$ est exceptionnel pour la paire $(X, \Delta +  D)$, c'est-à-dire $c_1 = c_2 = 0$ et $A = 0$. On peut donc supposer que $\tmop{codim}(W) > 1$.

La première étape consiste à rendre le centre $W$ maximal. Pour cela on pose $D_1=(1-\epsilon)D$ pour un rationnel  $\epsilon$ suffisamment petit. La paire $(X,\Delta + D_1)$ est klt. On choisit maintenant un $\mathbbm{Q}$-diviseur de Cartier ample $A_1\sim_{\mathbbm{Q}}H$ très général parmi ceux contenant $W$. Plus précisément, on suppose que ce diviseur $A_1$ ne contient pas de centre de singularités log-canoniques de $(X, \Delta +  D)$ contenant strictement $W$. Il existe un rationnel $a_1$ tel que $(X,\Delta + D_1 + a_1 A_1 )$ soit log-canonique mais non klt au point général de $W$, c'est-à-dire que $W$ soit un centre de singularités log-canoniques maximal pour la paire $(X,\Delta + D_1 + a_1 A_1 )$. On peut remarquer que $a_1 \rightarrow 0$ lorsque $\epsilon \rightarrow 0$.

On peut donc supposer dès maintenant que $W$ est un centre de singularités log-canoniques maximal pour la paire $(X, \Delta + D)$. En particulier, tout diviseur sur une log-résolution de $(X,\Delta +D)$ de discrépance $-1$ pour $(X,\Delta +D)$ et dominant $W$ a son image égale à $W$.

Soit $\mu : {X'} \rightarrow X$ une log-résolution de $(X,\Delta +D)$. On va construire un $\mathbbm{Q}$-diviseur ample sur ${X'}$ qui va nous permettre de garder un seul diviseur exceptionnel de discrépance $-1$ au-dessus de $W$. Toute la difficulté consiste à ce que ce $\mathbbm{Q}$-diviseur ample <<vienne d'en bas>>.

On note $a_i$ la discrépance de $E_i$ pour la paire $(X,\Delta +D)$, $d_i$ le coefficient de $E_i$ dans le $\mathbbm{Q}$-diviseur $\mu^* D$.

On commence par construire un diviseur ample sur $X'$.
Pour un diviseur $F$ dont le support est $\mu$-exceptionnel, on note $F^C$ la somme réduite des diviseurs $\mu$-exceptionnels non inclus dans le support de $F$. 
On note $D_{\tmop{exc}}$ le diviseur $\mu^* D - \mu^{-1}_* D$. Comme son nom l'indique $D_{\tmop{exc}}$ est supporté sur le lieu exceptionnel de $\mu$.

Pour tout rationnel $0 < \varepsilon \ll 1$, le $\mathbbm{Q}$-diviseur $\mu^* A_1 - \varepsilon (D_{\tmop{exc}} + D_{\tmop{exc}}^C)$ est ample. On note $B \sim_{\mathbbm{Q}} \mu^* A_1 - \varepsilon (D_{\tmop{exc}} + D_{\tmop{exc}}^C)$  un $\mathbbm{Q}$-diviseur effectif très général.

On va maintenant modifier légèrement notre diviseur ample $B$ sur $X$ de sorte que l'image directe dans $X$ du nouveau diviseur obtenu soit le diviseur recherché.
On choisit un diviseur $E_0$ sur $X'$ dominant $W$, de discrépance $-1$ pour la paire $(X,\Delta +D)$ et de coefficient $d_0$ minimal parmi de tels diviseurs :
$$d_0 = \min \{d_i \mid E_i \textrm{ soit un diviseur de discrépance } -1 \textrm{ dominant }W \textrm{ et de coefficient } d_i \textrm{ dans } \mu^* D\}.$$
On peut remarquer que la paire $(X,\Delta)$ étant klt, on a $a(E_0, X, \Delta) > -1$. On en déduit que $d_0 > 0$ puisque $a(E_0, X, \Delta + D) = -1$.
On pose $B' \sim_{\mathbbm{Q}} B - \varepsilon E_0$ de sorte que $B'$ soit ample, effectif, irréductible et que 
$B' +  \varepsilon E_0 + \varepsilon (D_{\tmop{exc}} + D_{\tmop{exc}}^C)$ soit à croisements normaux simples
: ceci est possible si $\varepsilon$ est choisi dès le départ suffisamment petit.

Par construction, on a $B' +  \varepsilon E_0 + \varepsilon (D_{\tmop{exc}} + D_{\tmop{exc}}^C) \sim  \mu^* A_1$. On pose $B_X=\mu_* (B' +  \varepsilon E_0 + \varepsilon (D_{\tmop{exc}} + D_{\tmop{exc}}^C))$. 
Ce $\mathbbm{Q}$-diviseur $B_X$ vérifie $B_X \sim A_1$ et est $\mathbbm{Q}$-Cartier. Comme $\mu^*B_X = B' + E'$, avec $E'$ exceptionnel et effectif, et que $\mu^*B_X \sim B' +  \varepsilon E_0 + \varepsilon (D_{\tmop{exc}} + D_{\tmop{exc}}^C) $, on en déduit que $\mu^*B_X = B' + \varepsilon (E_0 + D_{\tmop{exc}} + D_{\tmop{exc}}^C)$.

Pour tout rationnel $0 <\alpha \ll 1$, montrons qu'il existe un rationnel $\beta$ tel que $W$ soit un centre de singularités log-canoniques exceptionnel pour la paire $(X,\Delta +(1-\alpha)D + \beta B_X)$.
 On pose
\begin{eqnarray*}
d'_0 = 2 \varepsilon  & & \\
d'_i = \varepsilon & \textrm{si} & i\neq 0.
\end{eqnarray*} 

On a
$$K_{X'} + \mu^{-1}_* (\Delta +(1-\alpha)D + \beta B_X) = \mu^*(K_X + \Delta +(1-\alpha)D + \beta B_X) + \sum_{E_i \textrm{ }\mu\mathrm{-exc.}}
(a_i + \alpha d_i - \beta d'_i)E_i.$$
 
D'après l'équation précédente, si $0 <\alpha \ll 1$, $\beta = \frac{d_0}{2 \varepsilon}\alpha$ convient : en effet, on a
$a_0 + \alpha d_0 - \beta d'_0 = a_0 = -1$ et
$$a_i + \alpha d_i - \beta d'_i = a_i + \alpha d_i - \frac{\alpha}{2} d_0 \geqslant a_i + \frac{\alpha}{2} d_0 > a_i = -1$$
si $E_i$ est un diviseur dominant $W$ et vérifiant $i \neq 0$ et $a_i = -1$. Il reste le cas des diviseurs $E_j$ dominant $W$ et de discrépance $a_j > -1$. Mais ces diviseurs sont en nombre fini donc pour $\alpha$ suffisamment petit, on a $a_j + \alpha d_j - \beta d'_j > -1$.

Donc $E_0$ est le seul diviseur de discrépance $-1$ pour la paire $(X,\Delta +(1-\alpha)D + \beta B_X)$ dominant $W$. De plus $W$ est un centre de singularités log-canoniques minimal pour la paire $(X,\Delta +(1-\alpha)D + \beta B_X)$.
Cela conclut la démonstration en posant $c_1 = \beta$, $c_2 = \alpha$ et $A = \beta B_X$.

\end{proof}

\subsubsection{Sous-adjonction}

Une vaste généralisation du th\'eor\`eme d'adjonction pour les diviseurs lisses,
d\'emontr\'e par Kawamata dans {\cite{Kawamata1997}} et {\cite{Kawamata1998}},
permet de munir d'une structure de paire, d\'ependant du (log)-diviseur
canonique de $(X, D)$, les centres log-canoniques exceptionnels de la paire
$(X, D)$.

\begin{theorem}\label{sous-adjonction}
  {\dueto{{\cite{Kollar2007}}, theorem 8.6.1}}Soit $(X, D)$ une paire
  log-canonique et $W$ un centre de singularit\'es log-canoniques
  exceptionnel. Soit $H$ un diviseur ample et $\varepsilon > 0$ un rationnel.
  
  Alors $W$ est normal et il existe un $\mathbbm{Q}$-diviseur effectif $D_W$
  sur $W$ tel que :
  \begin{enumeratenumeric}
    \item $(W, D_W)$ soit une paire klt,
    
    \item $(K_X + D + \varepsilon H)_{|W} \sim_{\mathbbm{Q}} K_W + D_W$.
  \end{enumeratenumeric}
\end{theorem}

\subsection{Diviseurs de grand volume}\label{grandvolume}

On est maintenant à même de démontrer le théorème \ref{non-annulation}.

\begin{proof}\dueto{théorème \ref{non-annulation}}

  On note $x_0$ un point de $X$. D'apr\`es l'hypoth\`ese sur le volume de $L$,
  il existe un diviseur $D \sim_{\mathbbm{Q}} L$ de multiplicit\'e
  $\tmop{mult}_{x_0} D > 3.$ On en d\'eduit que le seuil log-canonique de $D$
  en $x_0$ v\'erifie $c = \tmop{lct} (D ; x_0) < 1$.
  
  Le $\mathbbm{Q}$-diviseur $L - cD$ est donc ample et d'apr\`es le
  th\'eor\`eme d'annulation de Nadel, la cohomologie sup\'erieure de
  $\mathcal{O}_X (K_X + L) \otimes \mathcal{J}(cD)$ est nulle. Si on note $V =
  Z (\mathcal{J}(cD))_{\textrm{red}}$ (on peut remarquer que $V \neq X$), on a alors une
  suite exacte
  \[ H^0 (X, \mathcal{O}_X (K_X + L) \otimes \mathcal{J}(cD)) \rightarrow H^0
     (X, \mathcal{O}_X (K_X + L)) \rightarrow H^0 (V, \mathcal{O}_V ((K_X +
     L)_{|V})) \rightarrow 0. \]
  L'existence d'une section globale non nulle de $\mathcal{O}_X (K_X + L)$
  d\'ecoule donc de l'existence d'une section globale non nulle de
  $\mathcal{O}_V ((K_X + L)_{|V})$. Si une composante irr\'eductible de $V$
  est de dimension nulle, c'est une composante connexe de $V$ et on a alors
  $H^0 (V, \mathcal{O}_V ((K_X + L)_{|V})) \neq 0$. On en d\'eduit donc que
  $H^0 (X, \mathcal{O}_X (K_X + L)) \neq 0$.
  
  Si toutes les composantes de $V$ sont de dimension strictement positive, on
  va l\'eg\`erement modifier le diviseur $D$ et son seuil log-canonique de
  sorte que l'une d'entre elles soit un centre de singularit\'es log-canonique
  exceptionnel.
  
  Soit $W$ un centre de singularit\'es log-canoniques minimal pour la paire $(X, cD)$. En
  appliquant le théorème \ref{perturbation} avec $\Delta = 0$, $H = L$, on obtient que $W$ est un centre log-canonique
  exceptionnel de $(X, c' D)$ pour un rationnel $c' < 1$. On peut
  donc appliquer le th\'eor\`eme de sous-adjonction (théorème \ref{sous-adjonction}) \`a la paire $(X, c' D)$.
  Pour tout rationnel $\varepsilon$, il existe un
  $\mathbbm{Q}$-diviseur effectif $\Delta_W$ sur $W$ v\'erifiant :
  \[ (K_X + c' D + \varepsilon D)_{|W} \sim_{\mathbbm{Q}} K_W + \Delta_W . \]
  La paire $(W, \Delta_W)$ étant klt et de dimension au plus $2$, on peut
  appliquer le th\'eor\`eme de non-annulation effective (\cite{Kawamata2000} theorem 3.1). 
  Le diviseur $(K_X + L)_{|W}$ est nef et
  \[ (K_X + L)_{|W} - (K_W + \Delta_W) \sim_{\mathbbm{Q}} (1 - c' - \varepsilon ) L_{|W}
     \]
     est ample si $\varepsilon$ est suffisamment petit.
  Il existe ainsi une section globale non nulle de $\mathcal{O}_W
  ((K_X + L)_{|W})$. De plus, $W$ \'etant une composante connexe de $V = Z
  (\mathcal{J}(c' D))_{\textrm{red}}$, on obtient donc une section globale non nulle de
  $\mathcal{O}_V ((K_X + L)_{|V})$. \`A l'aide de la suite exacte
  ci-dessus, on en d\'eduit l'existence d'une section globale non nulle de
  $\mathcal{O}_X (K_X + L)$.

\end{proof}

\subsection{Variétés presque de Fano}

Le point essentiel pour la d\'emonstration de la minoration des constantes de
Seshadri pour le diviseur fondamental de $- K_X$ (proposition
\ref{anticoind}) est l'existence d'une \'echelle pour la paire $(X, \Delta_X)$
({\cite{Ambro1999a}}, main theorem) :

\begin{theorem}
  {\dueto{Ambro}}\label{echelle}Soit $H$ un diviseur de Cartier gros et nef
  sur une vari\'et\'e projective normale $X$ de dimension $n$. Supposons qu'il
  existe une fronti\`ere $\Delta_X$ sur $X$ telle que
  \begin{enumerate}
    \item $(X, \Delta_X)$ soit klt,
    
    \item $- (K_X + \Delta_X) \equiv (n - c + 1) H$, pour un rationnel $c$
    v\'erifiant $n - c + 1 > 0$,
    
    \item $c < 4$.
  \end{enumerate}
  Alors $\dim |H| \geqslant n - 1$ et $|H|$ n'a pas de composante fixe. De
  plus, la paire $(X, \Delta_X + S)$ est purement log-terminale pour $S \in |H|$
  g\'en\'eral.
  
  En particulier, $(S, \Delta_{X|_S})$ est une paire presque de Fano de coïndice $c$
  si $n > c$, est Calabi-Yau si $c = n$ et est de type g\'en\'eral si $n < c <
  n + 1$.
\end{theorem}

\section{Constantes de Seshadri}
On exploite les résultats de non-annulation pour ramener la minoration des constantes de Seshadri au cas des surfaces lorsque c'est possible. La théorie de l'adjonction pour les variétés ayant des arêtes extrémales de grande longueur permet d'obtenir les minorations dans les cas restants.

\subsection{Variétés de dimension 3}
\subsubsection{\label{classificationdim3}Rappels lorsque le diviseur anticanonique est nef}

On s'int\'eressera dans cette sous-section aux vari\'et\'es lisses de dimension $3$
dont le diviseur anticanonique est nef et non trivial. On s'attardera plus
particuli\`erement sur les vari\'et\'es de dimension num\'erique \'egale \`a
$2$ et d'irr\'egularit\'e nulle.

\paragraph{Notations et rappels pr\'eliminaires}

La r\'eduction nef d\'efinie ci-dessous est une application presque
holomorphe. Commen\c{c}ons par d\'efinir cette derni\`ere notion.

\begin{definition}
  {\dueto{{\cite{Bauer2002}}, definition 2.3}}{\subindex{presque
  holomorphe}{application rationnelle $-$}}Soit $X$ et $Y$ des vari\'et\'es
  projectives normales $f : X \dashrightarrow Y$ une application rationnelle
  et $X^0$ l'ouvert maximal sur lequel $f$ est holomorphe. L'application $f$
  est dite presque holomorphe si certaines des fibres de la restriction
  $f_{|X^0}$ sont compactes.
\end{definition}

Cette d\'efinition \'equivaut \`a demander que l'image du lieu
d'ind\'etermination de l'application rationnelle $f$ ne domine pas $Y$.

\begin{definition}
  {\dueto{R\'eduction nef}}{\index{r\'eduction nef}}Soit $L$ un fibr\'e en
  droites nef sur une vari\'et\'e $X$ projective normale et de dimension $n$.
  La r\'eduction nef de $L$ est une application rationnelle, presque
  holomorphe et dominante $f : X \dashrightarrow B$ \`a fibres connexes telle
  que :
  \begin{enumerate}
    \item $L$ soit num\'eriquement trivial sur toutes les fibres compactes de
    dimension $\dim X - \dim B$,
    
    \item pour un point $x \in X$ en position g\'en\'erale et pour toute
    courbe irr\'eductible $C$ passant par $x$ v\'erifiant $\dim f (C) > 0$, on ait
    $L \cdot C > 0$.
  \end{enumerate}
\end{definition}

Dans {\cite{Bauer2002}}, les (nombreux) auteurs ont montr\'e l'existence et
l'unicit\'e d'une telle r\'eduction nef.

\begin{theorem}
  {\dueto{{\cite{Bauer2002}}}}Soit $L$ un fibr\'e en droites nef sur une
  vari\'et\'e normale projective $X$ de dimension n. La r\'eduction nef $f : X
  \dashrightarrow B$ de $L$ existe et est unique \`a \'equivalence
  birationnelle de $B$ pr\`es.
\end{theorem}

La dimension de la vari\'et\'e $B$ dans le th\'eor\`eme pr\'ec\'edent est
donc un invariant du fibré en droites nef $L$. \

\begin{definition}
  \label{defdimnef}La dimension nef{\index{dimension nef}} $n (L)$ d'un
  fibré en droites nef $L$ est la dimension de l'image de la r\'eduction
  nef de $L$. On d\'efinit la dimension
  nef d'un diviseur de Cartier nef $D$ en posant $n(D) = n(\mathcal{O}_X (D))$.
 La dimension nef d'une vari\'et\'e dont le
  diviseur anticanonique est nef est $n (X) = n (- K_X)$.
\end{definition}

On peut comparer la dimension nef d'un diviseur \`a sa dimension de Kodaira et
\`a sa dimension num\'erique.

\begin{proposition}
  \label{inegdim}{\dueto{Voir {\cite{Bauer2004}} definition-proposition
  1.4}}Soit $X$ une vari\'et\'e projective lisse et $D$ un diviseur nef. Si on
  note $\kappa (D)$ la dimension de Kodaira-Iitaka de $D$, $\nu (D)$ sa
  dimension num\'erique et $n (D)$ sa dimension nef alors
  \[ \kappa (D) \leqslant \nu (D) \leqslant n (D) . \]
\end{proposition}

Dans {\cite{Bauer2004}}, les auteurs \'etudient la r\'eduction nef du diviseur
anticanonique pour les vari\'et\'es de dimension $3$ \`a fibr\'e anticanonique
nef et non trivial. En comparant cette r\'eduction nef avec les contractions
des ar\^etes du c\^one de Mori ils en tirent plusieurs r\'esultats de
classification. Pour les vari\'et\'es de dimension $3$, la r\'eduction nef de
$- K_X$ est holomorphe.

\begin{theorem}
  \label{reducnefXdim3}{\dueto{{\cite{Bauer2004}}, theorem 2.1}}Soit $X$ une
  vari\'et\'e projective lisse de dimension $3$ dont le diviseur anticanonique
  $- K_X$ est nef. La r\'eduction nef de $- K_X$ est holomorphe et $- K_X$ est
  trivial sur chacune de ses fibres. De plus lorsque $X$ est rationnellement
  connexe et que $n (- K_X) = 1$ ou $n (- K_X) = 2$ alors $- K_X$ est
  semi-ample et la r\'eduction nef de $- K_X$ est la factorisation de Stein du
  morphisme induit par un grand multiple de $- K_X$.
\end{theorem}

\paragraph{Quelques \'el\'ements de classification}

On commence par rappeler un r\'esultat sur les vari\'et\'es d'irr\'egularit\'e
nulle.

\begin{proposition}
  \label{dimnumirreg0}Soit $X$ une
  vari\'et\'e projective lisse de dimension $3$ d'irr\'egularit\'e
  $q (X) = 0$ et de diviseur anticanonique $- K_X$ nef mais non
  num\'eriquement trivial (i.e $- K_X \nequiv 0$). Si $X$ n'est pas
  rationnellement connexe alors on a $n (- K_X) = 1$.
\end{proposition}

\begin{proof}
 Il s'agit de {\cite{Bauer2004}}, corollary 4.4 et {\cite{Bauer2004}}, corollary 3.2 combinés ensemble.
\end{proof}

On s'int\'eresse maintenant plus sp\'ecifiquement aux vari\'et\'es
rationnellement connexes, d'irr\'egularit\'e nulle et de dimension num\'erique
$\nu (X) = 2$. On a vu \`a la proposition \ref{inegdim} que dans ce cas $n
(X) = 2$ ou $n (X) = 3$.

Supposons dans un premier temps que $n (X) = 2$. La r\'eduction nef de $- K_X$
est un morphisme $f : X \rightarrow B$ vers une surface $B$ normale. On note
$\varphi : X \rightarrow Y$ une contraction \'el\'ementaire de Mori. Il en
existe au moins une puisque $- K_X \nequiv 0$. On peut remarquer que le
syst\`eme lin\'eaire $| - K_X |$ est non vide : $h^0 (- K_X) \geqslant 3$
(voir {\cite{Bauer2004}}, page 335). Commen\c{c}ons par le cas $\dim Y = 2$. Nous
n'aurons besoin que du cas d'un fibr\'e en coniques dont le discriminant est
vide.

\begin{theorem}
  \label{dimY=2}{\dueto{{\cite{Bauer2004}}, theorem 6.6}}Avec les notations
  ci-dessus, supposons que $\dim Y = 2$. Soit $\Delta$ le discriminant du
  fibr\'e en coniques $\varphi : X \rightarrow Y$ et supposons $\Delta =
  \varnothing$. Alors $\varphi$ est un fibr\'e en $\mathbbm{P}^1$ et $X
  =\mathbbm{P}(E)$ pour un fibr\'e vectoriel $E$ de rang $2$ sur $Y$. En
  particulier $- K_Y$ est nef. De plus on a soit $B =\mathbbm{P}^2$, soit $B
  =\mathbbm{P}^1 \times \mathbbm{P}^1$ ou soit $B$ est l'\'eclatement de
  $\mathbbm{P}^2$ en un point et
  \begin{enumerate}
    \item si $B =\mathbbm{P}^2$ alors $X \subset \mathbbm{P}^2 \times Y$ est
    donn\'e par $X \in |\mathcal{O}_{\mathbbm{P}^2} (1) \boxtimes
    \mathcal{O}_Y (- K_Y) |$ et $K_Y^2 > 0$,
    
    \item si $B =\mathbbm{P}^1 \times \mathbbm{P}^1$ ou si $B$ est
    l'\'eclatement de $\mathbbm{P}^2$ en un point, alors $Y$ est
    $\mathbbm{P}^2$ \'eclat\'e en $9$ points de sorte que $Y$ soit munie d'une
    fibration elliptique $g : Y \rightarrow \mathbbm{P}^1$ et que $E =
    g^{\ast} (\mathcal{O}_{\mathbbm{P}^1} (a) \oplus
    \mathcal{O}_{\mathbbm{P}^1})$ avec $a = 0$ ou $a = 1$.
  \end{enumerate}
\end{theorem}

On s'int\'eresse ensuite au cas $\dim Y = 3$.

\begin{proposition}
  \label{dimY=3}{\dueto{Voir {\cite{Bauer2004}}, proposition 6.7}}Toujours
  avec les m\^eme notations, supposons que $\dim Y = 3$. Notons $E$ le
  diviseur exceptionnel de l'\'eclatement $\varphi : X \rightarrow Y$ et
  supposons que $\dim \varphi (E) = 0$. Alors la vari\'et\'e $Y$ est terminale
  et le diviseur anticanonique $- K_Y$ est gros et nef.
\end{proposition}

Il reste maintenant \`a consid\'erer le cas des vari\'et\'es de dimension nef
\'egale \`a $3$. Les r\'esultats sur ces vari\'et\'es sont moins pr\'ecis mais
sont suffisants pour notre propos. On a notamment un r\'esultat sur la
structure du syst\`eme lin\'eaire associ\'e au diviseur anticanonique.

\begin{proposition}
  \label{n(X)=3}{\dueto{{\cite{Bauer2004}}, proposition 7.2}}Soit $X$ une
  vari\'et\'e projective de dimension $3$ rationnellement connexe de diviseur
  anticanonique $- K_X$ nef et v\'erifiant $n (- K_X) = 3$ et $\nu (- K_X) =
  2$. Alors le syst\`eme lin\'eaire $| - K_X |$ a une partie fixe non vide, c'est-à-dire une composante divisorielle de son lieu de base, que l'on note
  $A$. La partie mobile induit une fibration $f : X \rightarrow
  \mathbbm{P}^1$. Si $F$ est une fibre de $f$ alors $| - K_X | = A + |kF|$
  avec $k \geqslant 2$. De plus $A^3 = A^2 \cdot F = 0$.
\end{proposition}

\subsubsection{Preuves}

Pour les vari\'et\'es de dimension $3$, il est naturel d'essayer de se
ramener au cas des surfaces. On dispose en effet pour les surfaces lisses d'une minoration optimale des constantes de Seshadri des fibrés en droites amples due à Ein et Lazarsfeld \cite{Ein1993}. Les r\'esultats de cette section sont bas\'es sur
le lemme suivant :

\begin{lemma}
  \label{Lem}Soit $X$ une vari\'et\'e projective normale de dimension $3$ et
  $L$ un diviseur de Cartier gros et nef sur $X$. S'il existe un diviseur
  effectif de Cartier $D$ sur $X$ tel que $L - D$ soit pseudo-effectif, alors
  $\varepsilon (L ; x) \geqslant 1$ pour tout point $x$ en position tr\`es
  g\'en\'erale.
\end{lemma}

\begin{proof}
  Si $L - D$ est pseudo-effectif, $L - D$ a une intersection positive avec
  toute courbe mobile sur $X$.
  En particulier, il existe une composante irréductible $D_0$ du support de $D$ et un ensemble $V$,
  intersection d\'enombrable d'ouverts denses de $D_0$ tels que si $x
  \in V$, toute courbe $C$ passant par $x$ et non incluse dans $D_0$
  v\'erifie $(L - D) \cdot C \geqslant 0$, c'est-\`a-dire
  \[  L \cdot C \geqslant D \cdot C \geqslant \tmop{mult}_x C . \]
Cet ensemble $V$ est construit de la façon suivante : il existe une union dénombrable de sous-variétés strictes\footnote{On peut par exemple choisir le lieu de base restreint de $L-D$.} $Z_i$ de $X$ qui contient toute courbe ayant une intersection strictement négative avec $L-D$.
On fixe une composante irréductible $D_0$ du support de $D$ et on considère l'union notée $Z$ des sous-variétés $Z_i$ ne contenant pas $D_0$. On pose $V = D_0 \setminus Z$. L'ensemble $V$ est bien la restriction à $D_0$ d'une intersection dénombrable d'ouverts ayant une intersection non vide avec $D_0$. De plus toute courbe s'intersectant strictement négativement avec $L - D$ et non contenue dans $D_0$ est contenue dans $Z$ et ne peut donc passer par un point de $V$.

  Il ne reste donc plus qu'\`a consid\'erer les courbes incluses dans
  $\tmop{supp} D$. Si $x$ est un point lisse de $\tmop{supp} D$ et en position
  tr\`es g\'en\'erale dans $\tmop{supp} D$, alors d'après la
  proposition \ref{minsurfgrosnef} ci-dessous, on a $\varepsilon (L_{| \tmop{supp}
  D} ; x) \geqslant 1$ et donc pour toute courbe $C$ incluse dans $\tmop{supp}
  D$ et passant par $x$ on a
  \[ L \cdot C = L_{| \tmop{supp} D} \cdot C \geqslant \tmop{mult}_{x} C. \]
  On en d\'eduit que $\varepsilon (L ; x) \geqslant 1$ et par
  semi-continuit\'e inf\'erieure des constantes de Seshadri
  ({\cite{Lazarsfeld2004a}}, Exemple 5.1.11), la minoration vaut pour tout
  point $x$ en position tr\`es g\'en\'erale dans $X$.
\end{proof}

Le diviseur $D$ apparaissant dans le lemme pr\'ec\'edent n'ayant aucune
raison d'\^etre lisse, il est n\'ecessaire d'\'etendre le r\'esultat de Ein et
Lazarsfeld \`a des surfaces singuli\`eres \'eventuellement non normales et à des diviseurs gros et nef :

\begin{proposition}\label{minsurfgrosnef}
 Soit $S$ une surface projective éventuellement non normale, et $L$ un fibr\'e en droites gros et nef sur $S$. Pour tout point $x$ en position très générale, on a
  \[ \varepsilon (L ; x) \geqslant 1.\]
\end{proposition}

\begin{proof}

Quitte à résoudre les singularités de $S$ on peut supposer $S$ lisse, $D$ reste gros et nef sur une résolution des singularités de $S$. Le diviseur $L$ étant gros, si les constantes de Seshadri sont strictement inférieures à $1$ en tout point de $S$, il existe alors sur $S$ une famille non triviale
  de courbes $(C_t)_{t \in \Delta}$ param\'etr\'ees par un disque telles qu'en
  un point $x_t \in C_t$, chaque courbe v\'erifie $\tmop{mult}_{x_t} C_t > L
  \cdot C_t \geqslant 1$.
  On pose $m = \tmop{mult}_{\mu^{- 1} (x_t)}  \bar{C}_t$ pour $t$ g\'en\'eral.
  D'apr\`es {\cite{Ein1993}}, corollary 1.2, on a alors
  \[ C_t^2 \geqslant m (m - 1) . \]
  Par le th\'eor\`eme d'indice de Hodge, on a alors
  \[ m (m - 1) \leqslant (\mu^*L)^2 \bar{C}_t^2 \leqslant (\mu^* L \cdot \bar{C}_t)^2 \leqslant (m -
     1)^2 \]
  ce qui est une contradiction si $m > 1$.
\end{proof}

On aura besoin par la suite de la dénombrabilité de l'ensemble des points lisses où la constante de Seshadri est petite. La preuve est identique à celle de la proposition ci-dessus.

\begin{proposition}
  \label{surfsing}Soit $S$ une surface projective, \'eventuellement non
  normale, et $L$ un fibr\'e en droites ample sur $S$. L'ensemble des points
  lisses $x$ de $S$ tels que
  \[ \varepsilon (L ; x) < 1 \]
  est au plus d\'enombrable.
\end{proposition}

\begin{proof}
  Supposons qu'il existe sur $S$ une famille non triviale
  de courbes $(C_t)_{t \in \Delta}$ param\'etr\'ees par un disque telles qu'en
  un point $x_t \in C_t$, chaque courbe v\'erifie $\tmop{mult}_{x_t} C_t > L
  \cdot C_t \geqslant 1$. Supposons de plus que l'ensemble $\{x_t |t \in
  \Delta\}$ soit non contenu dans le lieu singulier de $S$. On choisit une r\'esolution des
  singularit\'es $\mu : \bar{S} \rightarrow S$ de sorte que $\mu$ soit
  un isomorphisme sur le lieu lisse de $S$. On consid\`ere la transform\'ee
  stricte $\bar{C}_t$ de la courbe $C_t$ par $\mu$. L'ensemble $\{x_t |t \in
  \Delta\}$ \'etant non inclus dans $\tmop{Sing}(S)$, pour $t$ g\'en\'eral, le point $x_t$
  est en dehors du lieu exceptionnel de $\mu$ et en ces points on a donc
  toujours
  \[ \mu^{\ast} L \cdot \bar{C}_t < \tmop{mult}_{_{\mu^{- 1} (x_t)}} \bar{C}_t
     . \]
  On pose $m = \tmop{mult}_{\mu^{- 1} (x_t)}  \bar{C}_t$ pour $t$ g\'en\'eral.
  D'apr\`es {\cite{Ein1993}}, corollary 1.2, on a alors
  \[ C_t^2 \geqslant m (m - 1) . \]
  Par le th\'eor\`eme d'indice de Hodge, on a alors
  \[ m (m - 1) \leqslant (\mu^*L)^2 \bar{C}_t^2 \leqslant (\mu^* L \cdot \bar{C}_t)^2 \leqslant (m -
     1)^2 \]
  ce qui est une contradiction si $m > 1$.
\end{proof}

\begin{proof}
  {\dueto{Th\'eor\`eme \ref{Th}, point 1}}

    On commence par traiter le cas g\'en\'eral, c'est-\`a-dire $(\nu
    (X), q (X)) \neq (2, 0)$.
    
    On sait d'apr\`es le th\'eor\`eme \ref{c2} que dans ce cas, la deuxi\`eme
    classe de Chern de $X$ est pseudo-effective. On en d\'eduit que $H^0 (X,
    L) \neq 0$ en utilisant le lemme \ref{secglob}. En appliquant le lemme
    \ref{Lem} au diviseur ample $L$ et au diviseur effectif $D \in |L|$, on
    obtient $\varepsilon (L ; x) \geqslant 1$.
    
    En particulier cela montre le r\'esultat pour les diviseurs amples des
    vari\'et\'es de Fano de dimension $3$.
    
    Il reste donc \`a minorer les constantes de Seshadri des fibr\'es en
    droites amples sur les vari\'et\'es v\'erifiant $(\nu (X), q (X)) = (2,
    0)$.
    
     La
     dimension num\'erique de $X$ valant $\nu (X) = 2$, on a $n (X) \geqslant
     2$ d'apr\`es l'in\'egalit\'e de la proposition \ref{inegdim}. On peut
     alors d\'eduire de la proposition \ref{dimnumirreg0} que la vari\'et\'e
     $X$ est rationnellement connexe.

    On distingue deux sous-cas, selon la dimension nef de $X$.
      Commen\c{c}ons par le cas $n (X) = 2$.

       Supposons que $L + K_X$ soit nef. Le syst\`eme lin\'eaire associ\'e au
      diviseur anticanonique \'etant non vide, pour un \'el\'ement $S \in | -
      K_X |$ et un point $x$ en position g\'en\'erale sur $\tmop{supp} S$, la
      constante de Seshadri de $L$ en $x$ est minor\'ee par $1$ : en effet le
      diviseur $L - S$ est nef par hypoth\`ese et il suffit alors d'appliquer
      le lemme \ref{Lem} \`a $L$ et $S$.
      
      Si $L + K_X$ n'est pas nef, il existe une courbe rationnelle $\Gamma$
      engendrant une ar\^ete $R =\mathbbm{R}_+ \cdot [\Gamma]$ de
      $\overline{\tmop{NE}} (X)_{K_X < 0}$ et v\'erifiant $(L + K_X) \cdot
      \Gamma < 0$. Puisque $L$ est ample, on a $L \cdot C \geqslant 1$ pour
      toute courbe $C$, donc l'ar\^ete $R$ est de longueur au moins $2$. On
      consid\`ere la contraction extr\'emale $\varphi : X \rightarrow Y$
      associ\'ee.
      
      Commen\c{c}ons par \'etudier le cas d'une contraction divisorielle. En
      dimension $3$, les contractions extr\'emales divisorielles de longueur
      au moins $2$ sont les \'eclatements d'un point lisse (\cite{Mori1982}, theorems
       3.3 et 3.4). 
      La vari\'et\'e
      $Y$ est alors terminale et de diviseur anticanonique gros et nef
      (proposition \ref{dimY=3}). Sa deuxi\`eme classe de Chern est donc
      pseudo-effective (th\'eor\`eme \ref{c2} - point 4) et on en d\'eduit
      que la seconde classe de Chern de $X$ est pseudo-effective (\cite{Xie2004}, lemme
       4.7). 

      On conclut comme dans le point $1$.
      
      Supposons maintenant que la contraction associ\'ee \`a l'ar\^ete
      $\mathbbm{R}_+ \cdot [\Gamma]$ soit une fibration. Par hypoth\`ese, $X$
      n'est pas Fano donc $Y$ ne peut \^etre un point. Si $Y$ est une courbe
      alors $\rho (X) = 2$ et donc, d'apr\`es \cite{Xie2004}, corollary 4.14,
      $c_2 (X)$ est pseudo-effective. On conclut encore une
      fois comme au point 1.
      
      Il reste finalement \`a traiter le cas o\`u $Y$ est une surface.
      L'ar\^ete $\mathbbm{R}_+ \cdot [\Gamma]$ \'etant de longueur $2$, le
      discriminant de la fibration est vide. La vari\'et\'e $X$ est donc un
      fibr\'e projectif sur $Y$ de la forme $\mathbbm{P}(E)$, pour un fibr\'e
      vectoriel de rang $2$ sur $Y$ (th\'eor\`eme \ref{dimY=2}). Plus
      pr\'ecis\'ement soit $X$ est une hypersurface de $\mathbbm{P}^2 \times
      Y$ donn\'ee par les z\'eros d'une section de
      $\mathcal{O}_{\mathbbm{P}^2} (1) \boxtimes K_Y^{- 1}$, soit $Y$
      admet une fibration elliptique $g : Y \rightarrow \mathbbm{P}^1$ et
      $X$ est de la forme $\mathbbm{P}(g^{\ast} (\mathcal{O}_{\mathbbm{P}^1}
      (b) \oplus \mathcal{O}_{\mathbbm{P}^1}))$ avec $b = 0$ ou $b = 1$.
      
      Dans le premier cas, en \'ecrivant $L \equiv_{\mathbbm{Q}} af^{\ast}
      \mathcal{O}(1) + \varphi^{\ast} D$ pour un diviseur de Cartier $D$ sur
      $Y$ et en consid\'erant l'intersection de $L$ avec une fibre de
      $\varphi$, on montre que $a = 1$ : en effet, l'image d'une fibre $l_y$
      de $\varphi$ par $f$ est une droite or $\varphi^{\ast} D \cdot l_y =
      0$ c'est-\`a-dire $af^{\ast} \mathcal{O}(1) \cdot l_y = L \cdot l_y$, et
      $L \cdot l_y \in \mathbbm{N}^{\ast}$ donc $a$ est un entier. De plus $X$
      n'est pas l'une des vari\'et\'es apparaissant au th\'eor\`eme
      \ref{valnef<gtr>n-1} : on en d\'eduit que $2 L + K_X$ est nef. Ceci
      combin\'e au fait que le diviseur $L + K_X$ n'est pas nef et que $K_X
      \cdot l_y = - 2$ (en calculant par exemple $K_X$ par adjonction), on en
      d\'eduit que $a = 1$. On a alors $\varphi^{\ast} D$ nef, puisque $2 L +
      K_X = 2 \varphi^{\ast} D$ est nef. Le diviseur $f^{\ast} \mathcal{O}(1)$
      \'etant effectif, on peut appliquer le lemme \ref{Lem} \`a $L$ et
      $f^{\ast} \mathcal{O}(1)$ et les constantes de Seshadri de $L$ sont donc
      minor\'ees par $1$ en tout point en position tr\`es g\'en\'erale.
      
      Dans le second cas, on se donne un diviseur $H$ dans la classe numérique de $\mathcal{O}_{\mathbbm{P}(E)}(1)$. On peut alors écrire $L \equiv aH + \varphi^{\ast} D$ pour
      un diviseur $D$ sur $Y$ et un entier $a>0$. Montrons que $D$ est ample. Puisque $E= g^{\ast} (\mathcal{O}_{\mathbbm{P}^1}
      (b) \oplus \mathcal{O}_{\mathbbm{P}^1})$, on a un plongement de $Y \simeq \mathbbm{P}(g^{\ast} \mathcal{O}_{\mathbbm{P}^1})$ dans $X$ et $H_{|Y} \equiv 0$. Mais $L_{|\mathbbm{P}(g^{\ast} \mathcal{O}_{\mathbbm{P}^1})}$ est ample donc $\varphi^{\ast} D_{|\mathbbm{P}(g^{\ast} \mathcal{O}_{\mathbbm{P}^1})}$ est ample, c'est-à-dire $D$ est ample.
      Soit $y \in Y$ un point tel que $\varepsilon
      (D ; y) \geqslant 1$ (qui existe d'après la minoration des constantes de Seshari sur les surfaces, voir \cite{Ein1993}) et $x \in X$ un point tel que $\varphi (x) = y$.
      Montrons que toute courbe $C$ passant par $x$ v\'erifie $L \cdot C \geqslant
      \tmop{mult}_x C$ : si $C$ est une fibre il n'y a rien \`a prouver puisque $C$ est lisse ;
      sinon, on a alors $\tmop{mult}_x C \leqslant \deg \varphi_{|C} \cdot
      \tmop{mult}_y \varphi (C)$ et donc $\tmop{mult}_x C \leqslant (\deg
      \varphi_{|C}) ( D \cdot \varphi (C)) \leqslant L \cdot C$ par la formule de
      projection.

      Il reste enfin le cas $n (X) = 3$.

      La seule vari\'et\'e projective lisse de dimension $3$ ayant une ar\^ete
      de longueur $4$ \'etant $\mathbbm{P}^3$ ({\cite{Cho2002}}, Corollary
      0.3), le diviseur $K_X + 3 L$ est nef.

      D'apr\`es la proposition \ref{n(X)=3} ci-dessus, la partie mobile du
      syst\`eme lin\'eaire anticanonique induit un morphisme $f : X
      \rightarrow \mathbbm{P}^1$. Si on note $F$ une fibre de ce morphisme
      alors $| - K_X | = A + |kF|$, o\`u $A$ est la composante fixe du
      syst\`eme et $k \geqslant 2$. Soit $x$ un point sur le support de $A$ et
      $F_x$ une fibre de $f$ passant par le point $x$. Toute courbe $C$
      passant par $x$ et non incluse dans $\tmop{supp} A \cup \tmop{supp} F_x$
      v\'erifie
      \[ - K_X \cdot C \geqslant A \cdot C + kF_x \cdot C \geqslant 3
         \tmop{mult}_x C \]
      et donc $L \cdot C \geqslant \tmop{mult}_x C$ puisque $3 L + K_X$ est
      nef.

      De plus, pour une fibre $F$ lisse intersectant le lieu lisse de
      $\tmop{supp} A$, l'intersection $\tmop{supp} A \cap F$ \'etant une
      courbe, il existe un point $x \in \tmop{supp} A \cap \tmop{supp} F$ tel
      que $\varepsilon (L_{| \tmop{supp} F} ; x) \geqslant 1$ et $\varepsilon
      (L_{| \tmop{supp} A} ; x) \geqslant 1$ : ceci d\'ecoule de la
      proposition \ref{surfsing}, les points lisses de $F$ et de $\tmop{supp}
      A$ ne v\'erifiant pas cette minoration pour les constantes de Seshadri
      de $L$ \'etant d\'enombrables.

      On a donc $L \cdot C \geqslant \tmop{mult}_x C$ pour toute courbe $C$
      passant par $x$. On en d\'eduit donc que $\varepsilon (L ; x) \geqslant
      1$.

\end{proof}

On finit par la minoration des constantes de Seshadri pour les fibr\'es
adjoints.

\begin{proof}
  {\dueto{Théorème \ref{Th}, point 2}}D'apr\`es {\cite{Kawamata2000}}, proposition
  4.1,
  \[ H^0 (X, \mathcal{O}_X (K_X + L)) \neq 0 \]
  et comme dans la preuve du th\'eor\`eme pr\'ec\'edent, on en d\'eduit que
  $\varepsilon (X, K_X + L  ; x) \geqslant 1$ pour tout point $x$ en position tr\`es
  g\'en\'erale.

\end{proof}

\begin{proof}{\dueto{Théorème \ref{Th}, point 3}}
  On se ram\`ene au cas d'une surface (voir proposition \ref{surfsing}) en trouvant un membre $S \in |K_X + L|$.
  
  Si $\sqrt[3]{L^3} > 3$, l'existence de $S$ est une application directe du
  th\'eor\`eme \ref{nonannulationdim3}.
  
  On suppose maintenant que $\sqrt[3]{L^3} = 3$, ce qui
  emp\^eche une application imm\'ediate du th\'eor\`eme
  \ref{nonannulationdim3}. On montre dans ce cas que soit
  \[ \varepsilon (X, K_X + L ; x) \geqslant \varepsilon (X, L ; x) \geqslant 1
  \]
  pour un point $x$ en position tr\`es g\'en\'erale, soit il existe un
  rationnel $0 < t < 1$ et un $\mathbbm{Q}$-diviseur $D
  \sim_{\mathbbm{Q}} tL$ tel que $\mathcal{J}(D) \neq \mathcal{O}_X$.
  
  Commen\c{c}ons par montrer que pour un point en position tr\`es g\'en\'erale
  $x$,
  \[ \varepsilon (X, K_X + L ; x) \geqslant \varepsilon (X, L ; x). \]
  En effet, puisque $K_X$ est pseudo-effectif, toute courbe passant par un point en position
  tr\`es g\'en\'erale $x$, v\'erifie $K_X \cdot C \geqslant 0$. On a donc
  $(K_X + L) \cdot C \geqslant L \cdot C$, ce qui implique l'in\'egalit\'e
  recherch\'ee.
  
  Soit $x \in X$ un point en position tr\`es g\'en\'erale, $f : X' \rightarrow
  X$ l'\'eclatement de $X$ en $x$, $E$ le diviseur exceptionnel. Supposons
  $\varepsilon (X, K_X + L ; x) < 1$. On note $c$ un rationnel v\'erifiant $1
  > c > \varepsilon (X, K_X + L ; x)$.
  Par hypothèse sur le volume de $L$, le diviseur $f^{\ast} L - (2 + c) E$
  est gros et d'après le choix de $c$, ce diviseur est non nef.
  On note $d = \tmop{lct} (\|f^{\ast} L - (2 + c) E\|)$.
  Si $d < 1$, il existe un rationnel $d' < 1$ et un $\mathbbm{Q}$-diviseur
  effectif $D' \sim_{\mathbbm{Q}} d' (f^{\ast} L - (2 + c) E)$ tel que
  $\mathcal{J}(D') \neq \mathcal{O}_{X'}$. D'apr\`es la formule de
  transformation birationnelle des id\'eaux multiplicateurs,
  $\mathcal{J}(f_{\ast} D') \neq \mathcal{O}_X$. Le diviseur $f_{\ast} D'$
  \'etant $\mathbbm{Q}$-lin\'eairement \'equivalent \`a $d' L$, en utilisant
  les m\^emes arguments qu'au cours de la preuve du th\'eor\`eme
  \ref{nonannulationdim3}, on en d\'eduit l'existence d'un diviseur effectif
  $S \in |K_X + L|$.
  
  Il reste maintenant le ``nouveau'' cas, lorsqu'il n'est pas possible de
  construire facilement un diviseur $D \sim_{\mathbbm{Q}} (1 - \alpha) L$
  v\'erifiant $\mathcal{J}(D) \neq \mathcal{O}_X$. On vient de voir que cela
  correspond au cas o\`u le seuil log-canonique de $\|f^{\ast} L - (2 + c)
  E\|$ v\'erifie $d = \tmop{lct} (\|f^{\ast} L - (2 + c) E\|) \geqslant 1$.
  
  L'id\'ee consiste alors \`a construire une courbe rationnelle dans $X$
  passant par $x$. La vari\'et\'e $X$ n'\'etant pas unir\'egl\'ee, on montre
  ainsi que le ``nouveau'' cas ne peut avoir lieu pour tout $x$. Pour ce
  faire, on va construire une fronti\`ere $\Delta$ telle que la paire $(X',
  \Delta)$ soit klt et $K_{X'} + \Delta$ soit non nef. On obtient ainsi par le
  th\'eor\`eme du c\^one une contraction, le lieu exceptionnel de cette contraction est
  unir\'egl\'e et on verra qu'il est localement isomorphe à son image sur $X$ et contient le point $x$.

Afin de mettre en lumière l'argument, traitons tout d'abord le cas où
$$d = \tmop{lct} (\|f^{\ast} L - (2 + c) E\|) > 1.$$
Le rationnel $c$ vérifie $c<1$, ce qui, par hypothèse sur $L$, implique que $2+c < \sqrt[3]{L^3}$. On peut donc choisir un $\mathbbm{Q}$-diviseur effectif vérifiant
  \[ \Delta \sim_{\mathbbm{Q}} f^{\ast} L - (2 + c) E.\]

Puisque $d>1$, on a alors $\mathcal{J}(\Delta) =\mathcal{O}_X$, autrement dit, la paire $(X', \Delta)$ est klt.
D'après notre choix de $c$ (on rappelle que $c> \varepsilon (X, K_X + L ; x)$), le diviseur $K_{X'} + \Delta = f^{\ast} (K_X+L) -c E$ n'est pas nef. 
Il existe donc une arête $(K_{X'} + \Delta)$-négative.  D'apr\`es le th\'eor\`eme du c\^one, cette arête peut
  \^etre contract\'ee et toute courbe dont la classe num\'erique appartient
  \`a cette arête est contract\'ee. Le lieu exceptionnel, noté $Z$, de cette contraction est
  unir\'egl\'e. De plus il est distinct de $E$.
  En effet, toute courbe contenue dans $E$ s'intersecte positivement avec $K_{X'}+ \Delta$. La restriction de $f : X' \rightarrow X$ à $Z$
  est donc localement un isomorphisme (en dehors de l'intersection de $Z$ avec $E$). En particulier l'image de $Z$ dans
  $X$ est uniréglée. Il reste à montrer qu'elle contient le point $x$. Puisque $K_X + L$ est ample, toute courbe dont la classe numérique est sur une arête $(K_{X'} + \Delta)$-négative a une intersection strictement positive avec $E$. Donc l'image de $Z$ dans $X$ contient $x$.

Le cas où $d=1$ est similaire. Pour un rationnel $0< \alpha \ll 1$, on construit
\[ \Delta \sim_{\mathbbm{Q}}  (1 - \alpha ) (f^{\ast} L - (2 + c) E).\]
Là encore, la paire $(X', \Delta)$ est klt. Il reste à montrer que $(K_{X'} + \Delta)$ n'est pas nef. On a
$$K_{X'} + \Delta = f^{\ast} (K_X+(1-\alpha) L) -(1 - \alpha) c E.$$
Pour toute courbe $C$ vérifiant $(f^{\ast} (K_X + L)-E) \cdot C <0$, il existe $0<\alpha \ll 1$ tel que
$$(1 - \alpha ) (K_{X'} + \Delta) \cdot C < 0.$$
Si $\alpha$ est suffisamment petit, il existe donc une arête $(K_{X'} + \Delta)$-négative. Le lieu exceptionnel de la contraction associée à cette arête est distinct de $E$ puisque toute courbe contenue dans $E$ s'intersecte positivement avec $K_X + \Delta$.
On conclut comme dans le cas $d>1$.

Puisque $X$ n'est pas uniréglé, il existe donc un point où la constante de Seshadri de $K_X + L$ est minorée par $1$. La borne est donc valable en tout point en position très générale.

\end{proof}

\subsection{Variétés presques Fano d'indice grand}
\subsubsection{Rappels sur la th\'eorie d'adjonction}

La d\'emonstration du th\'eor\`eme  \ref{coind-2} utilise la minoration des
constantes de Seshadri du diviseur anticanonique de la vari\'et\'e $X$ en
\'etudiant la valeur nef{\index{valeur nef}} du diviseur ample $A$. Cette
\'etude se fait classiquement dans le cadre des vari\'et\'es polaris\'ees $(X,
L)$, c'est-\`a-dire des couples form\'es d'une vari\'et\'e $X$ et d'un fibr\'e
en droites ample $L$.

\begin{definition}\label{valeurnef}
  Soit $A$ un diviseur ample sur une vari\'et\'e $X$ factorielle et \`a
  singularit\'es terminales. La valeur nef de $A$ est le rationnel $\tau (A) =
  \inf \{\tau \mid K_X + \tau A \text{ est nef} \}${\subindex{valeur
  nef}{diviseur}}.
\end{definition}

La d\'efinition ne d\'ependant que de la classe d'\'equivalence num\'erique de
$A$, pour un fibr\'e en droites ample{\subindex{valeur nef}{fibr\'e en
droites}} $L$, on d\'efinit la valeur nef de $L$ comme \'etant la valeur nef
de $c_1 (L)$. La valeur nef d'une vari\'et\'e polaris\'ee $(X, L)$ est la
valeur nef de la polarisation $L$.

On peut noter que le fait que la valeur nef de $A$ soit un rationnel est une
cons\'equence du th\'eor\`eme de rationalit\'e ({\cite{Debarre2001}}, theorem
7.34).

La th\'eorie de l'adjonction est l'\'etude des vari\'et\'es polaris\'ees de
grande valeur nef. On r\'esume dans le th\'eor\`eme \ref{valnef<gtr>n-1} la
classification des vari\'et\'es polaris\'ees $(X, L)$ dont la valeur nef est
strictement sup\'erieure \`a $\dim X - 1$. Auparavant, on donne quelques
d\'efinitions concernant les vari\'et\'es polaris\'ees que nous allons
rencontrer.

\begin{definition}
  {\index{c\^one g\'en\'eralis\'e}}{\dueto{C\^one g\'en\'eralis\'e}}Soit $L$
  un fibr\'e en droites tr\`es ample sur une vari\'et\'e projective $V$ de
  dimension $n$ et $E =\mathcal{O}_V^{\oplus N - n}$ o\`u $N \geqslant n$ est
  un entier. On note $\mathcal{C}=\mathbbm{P}(E \oplus L)$ et $\xi
  =\mathcal{O}_{\mathcal{C}} (1)$ le fibr\'e en droites tautologique sur
  $\mathcal{C}$. Soit $\varphi : \mathcal{C} \rightarrow \mathbbm{P}^m$ le
  morphisme d\'efini par les sections globales de $\xi$. On note
  $C_N (V, L)$ l'image de $\mathcal{C}$ par $\varphi$ et $\xi_L$ la
  restriction de $\mathcal{O}_{\mathbbm{P}^m} (1)$ \`a $C_N (V, L)$. La
  vari\'et\'e polaris\'ee $(C_N (V, L), \xi_L)$ est appel\'ee le c\^one
  g\'en\'eralis\'e de dimension $N$ de base $(V, L)$. On \'ecrira souvent de
  mani\`ere abusive $C_N (V, L)$ au lieu de $(C_N (V, L), \xi_L)$.
\end{definition}

\begin{definition}
  {\index{scroll}}{\dueto{Scroll}}Une vari\'et\'e polaris\'ee $(X, L)$
  $r$-Gorenstein (c'est-\`a-dire telle que $rK_X$ soit de Cartier pour un
  entier $r \geqslant 1$) de dimension $n$ est appel\'ee un scroll sur une
  vari\'et\'e $V$ de dimension $m$ s'il existe un morphisme surjectif \`a
  fibres connexes $p : X \rightarrow V$ et un diviseur de Cartier $A$ ample
  sur $V$ tel que
  \[ r (K_X + (n - m + 1) L) \sim p^{\ast} A. \]
\end{definition}

On peut noter que la fibre générale d'un scroll est un espace projectif (voir \cite{Beltrametti1993}).

\begin{definition}
  {\index{vari\'et\'e de del Pezzo}}{\dueto{Vari\'et\'es de del Pezzo}}Soit
  $X$ une vari\'et\'e normale $r$-Gorenstein de dimension $n$ dont le diviseur
  anticanonique est ample et $L$ un diviseur de Cartier ample sur $X$. On dit
  que $(X, L)$ est une vari\'et\'e de del Pezzo si $rK_X \sim - (n - 1) rL$.
\end{definition}

\begin{theorem}
  \label{valnef<gtr>n-1}{\dueto{{\cite{Beltrametti1995}}, proposition 7.2.2,
  theorem 7.2.3, theorem 7.2.4}}Soit $L$ un fibr\'e en droites ample sur une
  vari\'et\'e $X$ projective, normale, factorielle, \`a singularit\'es
  terminales et de dimension $n$. Supposons que la valeur nef de $L$ soit
  strictement sup\'erieure \`a $n - 1$. Alors la vari\'et\'e polaris\'ee $(X,
  L)$ est isomorphe \`a l'une des vari\'et\'es polaris\'ees suivantes :
  \begin{itemizeminus}
    \item $(\mathbbm{P}^n, \mathcal{O}(1))$ et $\tau (L) = n + 1$,
    
    \item $(Q, \mathcal{O}_Q (1))$, o\`u $Q$ est une hyperquadrique dans
    $\mathbbm{P}^{n + 1}$ et $\tau (L) = n$,
    
    \item $(\mathbbm{P}(E), \mathcal{O}_{\mathbbm{P}(E)} (1))$ pour un fibr\'e
    vectoriel $E$ de rang $n$ sur une courbe lisse et $\tau (L) = n$,
    
    \item $C_n (\mathbbm{P}^2, \mathcal{O}_{\mathbbm{P}^2} (2))$,
    c'est-\`a-dire un c\^one g\'en\'eralis\'e au-dessus de $(\mathbbm{P}^2,
    \mathcal{O}_{\mathbbm{P}^2} (2))$ et $\tau (L) = n - 1 / 2$.
  \end{itemizeminus}
\end{theorem}

Le r\'esultat pr\'ec\'edent apparaît aussi sous une forme synth\'etique dans
le livre de Beltrametti et Sommese {\cite{Beltrametti1995}}, table 7.1.

Bien s\^ur, plus la valeur nef de $L$ est petite, moins la paire $(X, L)$ est
contrainte. La premi\`ere valeur nef posant probl\`eme est le cas $\tau (L) =
n - 1$. En effet, si $X$ est l'\'eclatement d'une vari\'et\'e $Y$ projective
normale et factorielle en un point lisse, il existe un diviseur de Cartier
ample $L$ sur $X$ dont la valeur nef est $n - 1$. Il suffit de consid\'erer le
tir\'e en arri\`ere sur $X$ d'un diviseur de Cartier sur $Y$ suffisamment
ample et de lui soustraire une fois le diviseur exceptionnel. On souhaiterait
pouvoir aller plus loin dans l'analyse de la vari\'et\'e polaris\'ee $(X, L)$
et pour cela on consid\`ere le morphisme associ\'e \`a la valeur nef de $L$.
D'apr\`es le th\'eor\`eme d'absence de point de base (``base point freeness'') de
Kawamata-Shokurov, pour tout entier $k$ suffisamment divisible, le syst\`eme
lin\'eaire associ\'e \`a $k (K_X + \tau (L) L)$ est sans point base. Le
morphisme $\psi : X \rightarrow \mathbbm{P}^N$ associ\'e poss\`ede une
factorisation de Remmert-Stein $X \rightarrow Y \rightarrow \mathbbm{P}^N$ de
sorte que $X \rightarrow Y$ soit \`a fibres connexes, $Y$ soit normale et que
le morphisme $Y \rightarrow \mathbbm{P}^N$ soit fini.

\begin{definition}
  \label{morphismenef}On note $\varphi_L$ le morphisme $X \rightarrow Y$
  d\'efini ci-dessus. On l'appelle morphisme associ\'e \`a la valeur nef de
  $L$, qu'on abr\`egera parfois en morphisme nef de $L$.
\end{definition}

Le th\'eor\`eme suivant permet de distinguer les cas o\`u la valeur nef de $L$
est ``artificiellement'' gonfl\'ee \`a $n - 1$ par \'eclatement.

\begin{theorem}
  \label{valnef=n-1}{\dueto{{\cite{Beltrametti1995}} theorem 7.3.2}}Soit $X$
  une vari\'et\'e projective \`a singularit\'es terminales de dimension $n$ et
  $r$ un entier tel que $rK_X$ soit un diviseur entier. Soit $L$ un fibr\'e en
  droites ample sur $X$ et $\varphi_L$ le morphisme nef de $L$. Supposons que
  $\tau (L) \leqslant n - 1$. Alors $K_X + (n - 1)L$ est ample \`a moins que
  $\tau (L) = n - 1$ et qu'on soit dans l'un des cas suivants :
  \begin{enumerate}
    \item $- rK_X \sim r (n - 1) L$ et dans ce cas $(X, L)$ est une
    vari\'et\'e de del Pezzo,
    
    \item $(X, L)$ est une fibration en quadriques au-dessus d'une courbe
    lisse sous $\varphi_L$,
    
    \item $(X, L)$ est un scroll sous $\varphi_L$ au-dessus d'une surface normale,
    
    \item Le morphisme $\varphi_L : X \rightarrow Y$ est birationnel. Dans ce
    cas, si de plus $X$ est factorielle, alors $\varphi_L$ est la contraction
    simultan\'ee de diviseurs $E_i \simeq \mathbbm{P}^{n - 1}$ sur des points
    lisses distincts tels que $E_i \subset \tmop{reg} (X)$, $\mathcal{O}_{E_i}
    (E_i) \simeq \mathcal{O}_{\mathbbm{P}^{n - 1}} (- 1)$ et la restriction
    $L_{E_i}$ de $L$ \`a $E_i$ v\'erifie $L_{E_i} \simeq
    \mathcal{O}_{\mathbbm{P}^{n - 1}} (1)$. De plus, si on pose $L' =
    (\varphi_{L \ast} L)^{\ast \ast}$, alors $L'$ et $K_Y + (n - 1) L'$ sont
    amples et $K_X + (n - 1) L \sim \varphi^{\ast}_L (K_Y + (n - 1) L')$.
  \end{enumerate}
\end{theorem}

\begin{remark}
  Aux points 1, 2 et 3 du th\'eor\`eme pr\'ec\'edent, la restriction de $K_X +
  (n - 1) L$ \`a une fibre g\'en\'erale est num\'eriquement triviale. En
  particulier le diviseur $K_X + (n - 1) L$ n'est pas gros.
\end{remark}

On peut alors pousser plus loin la classification des vari\'et\'es
polaris\'ees de grande valeur nef en consid\'erant leur premi\`ere
r\'eduction.

\begin{definition}
  {\index{premi\`ere r\'eduction}}Soit une vari\'et\'e polaris\'ee $(X, L)$
  factorielle et \`a singularit\'es terminales. Si $K_X + (n - 1) L$ est gros
  et nef, on dit que $(X, L)$ admet une premi\`ere r\'eduction. Cette dernière est
une variété polarisée $(Y,L')$ telle que :
\begin{enumerate}
 \item si $K_X + (n -
  1) L$ est ample alors la premi\`ere r\'eduction de $(X, L)$ est $(Y,L')=(X,
  L)$,
\item si $K_X + (n - 1) L$ n'est pas ample, la premi\`ere r\'eduction de
  $(X, L)$ est la vari\'et\'e polaris\'ee $(Y, L')$ introduite au point 4 du
  pr\'ec\'edent th\'eor\`eme.
\end{enumerate}

\end{definition}

Apr\`es premi\`ere r\'eduction, on peut continuer la classification des
vari\'et\'es polaris\'ees. Le r\'esultat suivant est une version condens\'ee
de {\cite{Beltrametti1995}}, theorem 7.3.4.

\begin{theorem}
  {\dueto{Voir {\cite{Beltrametti1995}}, theorem 7.3.4}}\label{1erereduc}Soit
  $(X, L)$ une vari\'et\'e polaris\'ee normale factorielle et \`a
  singularit\'es terminales de dimension $n \geqslant 3$ admettant une
  premi\`ere r\'eduction $(Y, L')$. On note $\varphi_{L'}$ le morphisme
  associ\'e \`a la valeur nef de $L'$. Supposons que la valeur nef de $L'$
  v\'erifie $n - 2 < \tau (L') < n - 1$. Alors
  \begin{enumerate}
    \item soit $n = 4$, $\tau (L') = 5 / 2$, $(Y, L') \cong (\mathbbm{P}^4,
    \mathcal{O}_{\mathbbm{P}^4} (2))$,
    
    \item soit $n = 3$.
  \end{enumerate}
\end{theorem}

\subsubsection{Preuves}

Pour les vari\'et\'es presque de Fano de petit coïndice, l'existence d'une
\'echelle 
 pour le diviseur fondamental
permet de ramener la minoration des constantes de Seshadri au cas des petites
dimensions

\begin{proposition}
  \label{anticoind}
  Soit $(X, \Delta_X)$ une paire presque de Fano klt
  de dimension $n$, $H$ un diviseur de Cartier gros et nef sur $X$
  v\'erifiant $- (K_X + \Delta_X) \equiv (n - c + 1) H$ pour un rationnel $c < 4$.
  Alors pour un point $x$ en position tr\`es g\'en\'erale, on a
  \[ \varepsilon (X, H ; x) \geqslant 1. \]

\end{proposition}

Si de plus la vari\'et\'e $X$ est de Fano, factorielle et \`a singularit\'es
terminales, on peut 
obtenir une minoration des constantes de Seshadri de n'importe quel diviseur ample
$A$ sur $X$. On peut toutefois noter que dans le cas o\`u $X$ est lisse, de
Fano de dimension au moins $7$ et de coïndice au moins $4$, le nombre de
Picard de $X$ est $1$ ({\cite{Wisniewski1990}}, theorem A) : dans ce cas le théorème \ref{coind-2} n'apporte rien par rapport à la proposition \ref{anticoind}.

Lorsque $\tmop{Pic} (X)
  =\mathbbm{Z}$, une minoration des constantes de Seshadri du diviseur fondamental de $- K_X$
induit par homog\'en\'eit\'e une minoration pour tout diviseur ample.
Bien s\^ur, si $X$ est presque de Fano mais n'est pas de Fano, le nombre de
Picard de $X$ est strictement sup\'erieur \`a $1$, quelle que soit la
dimension de $X$. Dans le cas singulier, le nombre de Picard de $X$ n'est plus
aussi contraint par l'indice de la vari\'et\'e et il n'est pas clair \`a
l'heure actuelle{\footnote{Merci \`a Cinzia Casagrande d'avoir bien voulu
m'apporter des pr\'ecisions sur ce sujet.}} si un analogue de la conjecture de
Mukai {\cite{Mukai1988}} existe dans ce cas. Sur ce sujet, on pourra consulter
l'article {\cite{Casagrande2006}} de Casagrande, Jahnke et Radloff qui traite
le cas de la dimension $3$. Jahnke et Peternell donnent une
classification{\footnote{Nous n'utiliserons pas cette derni\`ere.}} des
vari\'et\'es presque de del Pezzo dans {\cite{Jahnke2006}}.

\paragraph{Le cas du diviseur anticanonique}

On peut maintenant prouver la proposition \ref{anticoind}.

\begin{proof}
  {\dueto{Proposition \ref{anticoind}}}
  
    L'id\'ee principale consiste à appliquer successivement le
  th\'eor\`eme \ref{echelle} \`a $X$ puis aux diviseurs ainsi obtenus et à se
  ram\`ener \`a minorer les constantes de Seshadri de $L$ sur une surface.
  
    On proc\`ede par r\'ecurrence. Supposons que pour toute paire $(Y, \Delta_Y)$ de
  dimension $n - 1$ v\'erifiant les hypoth\`eses de la proposition
  \ref{anticoind} on ait $\varepsilon (Y, H_Y ; x) \geqslant 1$, pour un point
  $x$ en position tr\`es g\'en\'erale. Soit $(X, \Delta_X)$ une paire de dimension
  $n \geqslant 4$ v\'erifiant les hypoth\`eses. D'apr\`es le th\'eor\`eme
  $\ref{echelle}$, il existe un \'el\'ement $S \in |H|$ irr\'eductible et
  r\'eduit tel que $(S, \Delta_{X|S})$ soit une paire presque de Fano. D'apr\`es
  notre hypoth\`ese de r\'ecurrence, les constantes de Seshadri de $H_{|S}$
  sont minor\'ees par $1$ pour tout point $y \in S$ en position tr\`es
  g\'en\'erale. Toute courbe $C$ incluse dans $S$ et passant par $y$ v\'erifie
  donc $H \cdot C = H_{|S} \cdot C \geqslant \tmop{mult}_y C$. De plus, le
  diviseur $S$ \'etant effectif, toute courbe $C$ non incluse dans $S$ et
  passant par un point $y \in \tmop{supp} S$ v\'erifie $H \cdot C = S \cdot C
  \geqslant \tmop{mult}_y C$. On en d\'eduit que $\varepsilon (X, H ; y)
  \geqslant 1$ et donc que $\varepsilon (X, H ; x) \geqslant 1$ pour tout
  point $x$ de $X$ en position tr\`es g\'en\'erale.
  
  Si $(X, \Delta_X)$ est de dimension $3$, il existe d'apr\`es le th\'eor\`eme
  \ref{echelle} un diviseur $S \in |H|$ irr\'eductible et r\'eduit. D'apr\`es
  le lemme \ref{Lem}, appliqu\'e \`a $H$ et $S$, les constantes de Seshadri de
  $H$ sont minor\'ees par $1$ en tout point $x$ en position g\'en\'erale.

\end{proof}

\paragraph{Vari\'et\'es presque Fano d'indice au moins $n - 2$}

Les cas difficiles pour la minoration des constantes de Seshadri du diviseur
ample $A$ sont ceux o\`u la valeur nef de $A$ est strictement sup\'erieure \`a
$n - 2$ : en effet, dans ce cas on ne peut utiliser de comparaison avec les constantes de Seshadri du diviseur fondamental du diviseur anticanonique.
On a vu cependant dans la sous-section pr\'ec\'edente que les vari\'et\'es
polaris\'ees de valeur nef strictement sup\'erieure \`a $n - 2$ sont
classifi\'ees, si l'on excepte celles admettant une premi\`ere r\'eduction.
Les deux lemmes suivants \'etablissent la minoration souhait\'ee des
constantes de Seshadri pour les vari\'et\'es intervenant dans cette
classification.

\begin{lemma}
  \label{minvalnef<gtr>n-1}Soit $(X, L)$ une des vari\'et\'es polaris\'ees du
  th\'eor\`eme \ref{valnef<gtr>n-1}. Supposons de plus que $X$ soit
  rationnellement connexe. Alors les constantes de Seshadri de $L$ sont
  minor\'ees par $1$ pour tout point $x \in X$ en position g\'en\'erale.
\end{lemma}

\begin{proof}

  Dans le cas de $\mathbbm{P}^n$ et de $Q_n$ et du c\^one g\'en\'eralis\'e sur
  $\mathbbm{P}^2$, le fibr\'e en droites $L$ est tr\`es ample donc les
  constantes de Seshadri de $L$ sont minor\'ees par $1$ en tout point de $X$.
  Dans le cas d'un fibr\'e projectif au-dessus d'une courbe lisse, on peut
  noter que la courbe est $\mathbbm{P}^1$, l'image d'une vari\'et\'e
  rationnellement connexe \'etant rationnellement connexe. Les constantes de
  Seshadri de $L$ sont alors minor\'ees par $1$ en tout point de $X$ en effet,
   $X$ est torique donc $L$ est sans point-base d'où la minoration souhaitée (\cite{Lazarsfeld2004a} example 5.1.18).

\end{proof}

\begin{lemma}
  \label{minvalnef=n-1}Soit $(X, L)$ une des vari\'et\'es apparaissant au
  point 2 ou au point 3 du th\'eor\`eme \ref{valnef=n-1}. Si l'on suppose de
  plus que $X$ est rationnellement connexe alors la constante de Seshadri de
  $L$ en un point $x$ en position très g\'en\'erale est minor\'ee par $1$.
\end{lemma}

\begin{proof}

  On suppose dans un premier temps que $(X, L)$ est une fibration en
  quadriques au-dessus d'une courbe lisse. La vari\'et\'e $X$ \'etant
  rationnellement connexe, cette courbe est $\mathbbm{P}^1$. La courbe \'etant
  lisse, la fibration est plate et le faisceau $\varphi_{L \ast} L$ est
  localement libre de rang $n + 2$. On a donc un plongement de $X$ dans
  $\mathbbm{P}(\varphi_{L \ast} L)$ de sorte que $L$ soit la restriction du
  fibr\'e $\mathcal{O}_{\mathbbm{P}(\varphi_{L \ast} L)} (1)$ \`a $X$ (voir
  aussi {\cite{Andreatta1993}}). Les constantes de Seshadri de
  $\mathcal{O}_{\mathbbm{P}(\varphi_{L \ast} L)} (1)$ \'etant minor\'ees par
  $1$ en tout point de $\mathbbm{P}(\varphi_{L \ast} L)$, il en va de m\^eme
  pour celles de $L$.
  
  Supposons maintenant que $(X, L)$ est un scroll au-dessus d'une surface
  normale $S$. On note $p : X \rightarrow S$ la projection sur $S$. La fibre
  g\'en\'erale \'etant un espace projectif, pour un point $x \in X$ en
  position tr\`es g\'en\'erale et une courbe $C$ passant par $x$ et incluse
  dans la fibre contenant $x$ on a $L \cdot C \geqslant \tmop{mult}_x C$.
  
  La vari\'et\'e $(X, L)$ \'etant un scroll au-dessus d'une surface, il existe
  un fibr\'e en droites ample $A$ sur $S$ v\'erifiant $K_X + (n - 1) L =
  p^{\ast} A$. Soit maintenant une courbe $C$ passant par $x$ non incluse dans
  une fibre de la projection $p$. Le point $x$ \'etant en position tr\`es
  g\'en\'erale dans $X$, il en est de m\^eme de son image dans $S$. On note
  $C'$ l'image de $C$ dans $S$ et $\nu : \overline{C'} \rightarrow C'$ et
  $\eta : \bar{C} \rightarrow C$ leurs normalisations respectives. Le
  morphisme $f = p_{|C} \circ \eta : \bar{C} \rightarrow C'$ se factorise en
  $\bar{C} \rightarrow \overline{C'} \rightarrow C'$. On note $d$ le degr\'e
  du morphisme $\bar{C} \rightarrow \overline{C'}$. De part la position
  g\'en\'erale de $p (x)$ dans $S$, on a $A \cdot C' \geqslant \tmop{mult}_{p
  (x)} C'$. De plus on a l'in\'egalit\'e suivante : $\tmop{mult}_x C \leqslant
  d \tmop{mult}_{p (x)} C'$. On a alors $p^{\ast} A \cdot C \geqslant d
  \tmop{mult}_x C' \geqslant \tmop{mult}_x C.$ De plus, d'apr\`es la
  proposition \ref{anticoind}, on a $H \cdot C \geqslant \tmop{mult}_x C$. On
  a donc
  \[ (n - 1) L \cdot C \geqslant - K_X \cdot C + p^{\ast} A \cdot C \geqslant
     rH \cdot C + p^{\ast} A \cdot C \geqslant (n - 2) \tmop{mult}_x C +
     \tmop{mult}_x C \]
  et on obtient finalement l'in\'egalit\'e souhait\'ee : $L \cdot C \geqslant
  \tmop{mult}_x C$. On a donc minor\'e les constantes de Seshadri de $L$ par
  $1$ pour tout point $x$ en position tr\`es g\'en\'erale de $X$.

\end{proof}

Apr\`es avoir trait\'e tous ces cas particuliers, on peut enfin d\'emontrer le
th\'eor\`eme principal :

\begin{proof}
  {\dueto{Th\'eor\`eme \ref{coind-2}}}
  On note $\tau (L) = \inf \{\tau \mid K_X +
  \tau L \text{ est nef} \}$ la valeur nef de $(X, L)$. On distingue les cas
  $\tau (L) \leqslant n - 2$ et $\tau (L) > n - 2$.
  
  Si $\tau (L) \leqslant n - 2$ alors pour tout point $x \in X$ en position
  tr\`es g\'en\'erale, $\varepsilon (X, L ; x) \geqslant 1$. Soit en effet $C$
  une courbe passant par $x$. Puisque $K_X + (n - 2) L$ est nef, on a
  \[ (n - 2) L \cdot C \geqslant - K_X \cdot C. \]
  Or $\varepsilon (X, - K_X ; x) \geqslant r_X \geqslant n - 2$ d'apr\`es la
  proposition \ref{anticoind} et l'hypoth\`ese sur l'indice de $X$. On en
  d\'eduit que $L \cdot C \geqslant \tmop{mult}_x C$, c'est-\`a-dire que
  $\varepsilon (X, L ; x) \geqslant 1$.
  
  Il reste \`a traiter le cas $\tau (L) > n - 2$. D'apr\`es les lemmes
  \ref{minvalnef<gtr>n-1} et \ref{minvalnef=n-1}, on peut supposer que $(X,
  L)$ admet une premi\`ere r\'eduction. Via cette r\'eduction, on obtient une
  vari\'et\'e polaris\'ee $(Y, L')$. D'apr\`es le point 4 du th\'eor\`eme
  \ref{valnef=n-1}, la valeur nef de $L'$ est strictement inf\'erieure \`a $n
  - 1$.
  
  Supposons que $n - 2 < \tau (L') < n - 1$. On peut alors utiliser le
  th\'eor\`eme \ref{1erereduc}.
  
  Les constantes de Seshadri de $\mathcal{O}_{\mathbbm{P}^4} (2)$ sont
  minor\'ees par $1$ en tout point. De plus, un \'eclatement de
  $\mathbbm{P}^4$ ne peut \^etre une vari\'et\'e de Fano d'indice sup\'erieur
  \`a $2$ : en effet si $X \rightarrow \mathbbm{P}^4$ est un \'eclatement de
  $\mathbbm{P}^4$ le long d'une sous-vari\'et\'e lisse, on a $\tmop{Pic} (X) =
  \tmop{Pic} (\mathbbm{P}^4) \oplus \mathbbm{Z} \cdot E$, o\`u $E$ est le
  diviseur exceptionnel. Donc si la premi\`ere r\'eduction de $(X, L)$ est
  $(\mathbbm{P}^4, \mathcal{O}_{\mathbbm{P}^4} (2))$, la premi\`ere
  r\'eduction est un isomorphisme : en effet, la premi\`ere r\'eduction est
  l'\'eclatement d'un nombre fini de points lisses. Dans ce cas, $L \simeq
  \mathcal{O}_{\mathbbm{P}^4} (2)$ et on a donc minor\'e les constantes de
  Seshadri de $L$ par $1$ en tout point de $X$.
  
  Si $\dim Y = 3$, alors la dimension de $X$ est ausi \'egale \`a $3$, le
  morphisme $f : X \rightarrow Y$ \'etant birationnel. La vari\'et\'e $X$ est
  donc une vari\'et\'e presque de Fano terminale. D'apr\`es le th\'eor\`eme
  \ref{c2} et le lemme \ref{secglob}, les constantes de Seshadri de $L$ sont
  minor\'ees par $1$ pour tout point $x$ en position tr\`es g\'en\'erale dans
  $X$.
  
  Si $\tau (L') \leqslant n - 2$ alors $K_Y + (n - 2) L'$ est nef. D'apr\`es
  le point 4 du th\'eor\`eme \ref{valnef=n-1}, on a $K_X = f^{\ast} K_Y + (n -
  1) E$, o\`u $E$ est le diviseur exceptionnel de $f : X \rightarrow Y$. De
  plus le diviseur $L$ est \'egal \`a $L = f^{\ast} L' - E$. On peut donc
  \'ecrire
  \[ K_X + (n - 2) L = f^{\ast} (K_Y + (n - 2) L') + E. \]
  Pour un point $x \in X$ en dehors du support de $E$ et toute courbe $C$
  passant par $x$, on alors
  \[ (K_X + (n - 2) L) \cdot C = f^{\ast} (K_Y + (n - 2) L' + E) \cdot C
     \geqslant 0. \]
  On peut donc conclure de la m\^eme mani\`ere que dans le cas $\tau (L)
  \leqslant n - 2$.

\end{proof}

\paragraph{Minoration pour le diviseur anticanonique des vari\'et\'es de Fano de
 dimension $4$}

\begin{proof}
  {\dueto{Th\'eor\`eme \ref{fano4}}}
  
  On commence par montrer que les constantes
  de Seshadri de $- K_X$ sont minor\'ees par $1$ en tout point en position
  tr\`es g\'en\'erale. D'apr\`es {\cite{Kawamata2000}}, theorem 5.2, il existe
  un diviseur effectif $Y \in | - K_X |$ tel que la paire $(X, Y)$ soit
  purement log terminale. En particulier, le diviseur $Y$ est irr\'eductible,
  normal et r\'eduit. De plus $Y$ est Gorenstein, ne poss\`ede que des
  singularit\'es canoniques et son diviseur canonique est lin\'eairement
  \'equivalent \`a $0$. D'apr\`es {\cite{Kawamata2000}}, proposition 4.1, le
  fibr\'e en droites $\mathcal{O}_Y (Y_{|Y})$ poss\`ede une section globale.
  On note $D$ le diviseur effectif associ\'e. Le diviseur $- K_{X|Y} - D$
  \'etant num\'eriquement trivial, on peut appliquer le lemme \ref{Lem} \`a $-
  K_{X|Y}$. Soit $x$ un point lisse de $Y$ pour lequel $\varepsilon (Y, Y_{|Y}
  ; x) \geqslant 1$. Toute courbe $C \subset X$ passant par $x$ v\'erifie $-
  K_X \cdot C \geqslant \tmop{mult}_x C$ : si $C$ est incluse dans $Y$ cela
  d\'ecoule de la minoration de la constante de Seshadri de $Y_{|Y}$ en $x$.
  Dans le cas contraire on a $- K_X \cdot C = Y \cdot C \geqslant
  \tmop{mult}_x C$. On a donc $\varepsilon (X, - K_X ; x) \geqslant 1$ et par
  semi-continuit\'e des constantes de Seshadri, les constantes de Seshadri de
  $- K_X$ sont minor\'ees par $1$ en tout point en position g\'en\'erale.

Cela conclut le cas où $r=1$, c'est-à-dire $-K_X=H$. Si $r>1$, la variété $X$
satisfait les hypothèses du théorème \ref{coind-2}.

\end{proof}

\newcommand{\etalchar}[1]{$^{#1}$}

\bigskip
 \noindent A.B. \emph{e-mail: broustet@ujf-grenoble.fr}
 
 \smallskip
 
 \noindent\emph{Institut Fourier, UMR 5582, UFR de Math\'ematiques\\Universit\'e Joseph Fourier,
 Grenoble 1\\BP 74\\38402 Saint Martin d'H\`eres, FRANCE}
\end{document}